\documentclass[11pt]{article}
\usepackage{graphicx}

\usepackage{amsmath}
\usepackage{amsfonts}
\usepackage{amssymb}
\usepackage{epsfig}

\usepackage{amsmath}
\usepackage{amsfonts}
\usepackage{amssymb}
\usepackage{dsfont}
\usepackage{mathrsfs}
\usepackage{bbm}
 \usepackage{hyperref}
\usepackage{amsmath}
\usepackage{amsfonts}
\usepackage{amssymb}
\usepackage{multirow}
\usepackage{mathrsfs}
\usepackage[dvipsnames,usenames]{color}

\usepackage{subcaption}
\usepackage{caption}
\usepackage{svg}

\usepackage{natbib}
\setcitestyle{numbers,square}

\renewcommand{\Box}{\framebox{\rule{0.3em}{0.0em}}}

\newtheorem{theorem}{Theorem}[section]
\newtheorem{theorem*}{Theorem}[subsubsection]
\newtheorem{lemma}{Lemma}[section]
\newtheorem{proposition}{Proposition}[section]
\newtheorem{example}{Example}[section]
\newtheorem{remark}{Remark}[section]
\newtheorem{definition}{Definition}[section]

\newtheorem{assumption}{Assumption}[section]

\renewcommand\thesubsubsection{\Alph{section}}

\newcommand{\setd}{{ d \kern -.15em l}}
\newcommand{\hatsetd}{ d \hat{\kern -.15em l }}
\newcommand{\dd}{\mathsf {d\kern -0.07em l}} 

\newcommand{\bgeqn}{\begin{eqnarray}}
\newcommand{\edeqn}{\end{eqnarray}}
\newcommand{\bgeq}{\begin{eqnarray*}}
\newcommand{\edeq}{\end{eqnarray*}}
\newcommand{\bec}{\begin{center}}
\newcommand{\enc}{\end{center}}
\newcommand{\R}{{\rm I\!R}}

\newcommand{\inmat}[1]{\mbox{\rm {#1}}}

\newcommand{\B}{{\cal B}}

\newcommand{\be}{\begin{equation}}
\newcommand{\ee}{\end{equation}}

\renewcommand{\Box}{\hfill \rule{2.3mm}{2.3mm}}

\title{
Generalized Bayesian Nash Equilibrium with Continuous Type and Action Spaces
}

\author{
Yuan Tao\thanks{Department of Systems Engineering and Engineering Management, The Chinese University of Hong Kong, Hong Kong, China. Email: ytao@se.cuhk.edu.hk
}
\,\,and\,\,
Huifu Xu\footnote{Corresponding author, Department of Systems Engineering and Engineering Management, The Chinese University of Hong Kong, Hong Kong, China. Email: hfxu@se.cuhk.edu.hk.}
}

\begin{document}

\maketitle

\begin{abstract}
{\color{black}
A Bayesian game is a strategic decision-making model where
each player's type parameter characterizing its own objective is private information: each player knows its own type but not its rivals' types, and a Bayesian Nash equilibrium (BNE) is an outcome of this game where each player makes a strategic optimal decision according to its own type under the Nash conjecture. In this paper, we advance the literature by considering a generalized Bayesian game where each player's action space depends on its own type parameter and its rivals' actions. This reflects the fact that in practical applications, a player's feasible action is often related to its own type (e.g.~marginal cost) and the rivals' actions (e.g.~common resource constraints in a competitive market). Under some moderate conditions, we demonstrate the existence of a continuous generalized Bayesian Nash equilibrium (GBNE) and the uniqueness of such an equilibrium when each player's action space is only dependent on its type. In the case that each player's action space also depends on rivals' actions, we give a simple example to show that   the uniqueness of GBNE is not guaranteed under standard monotone conditions. To compute an approximate GBNE, we restrict each player's response function to the space of polynomial functions of its type parameter and subsequently convert the GBNE model to a stochastic generalized Nash equilibrium (SGNE) model. To justify the approximation, we discuss the convergence of the approximation scheme. Some preliminary numerical test results show that the approximation scheme works well.
}
\end{abstract}


\textbf{MSCcodes.} {91A27, 91A06, 91A10, 91A15, 90C31}

\section{Introduction}

{\color{black}
A Bayesian game is a strategic interactive 
decision-making model where each player's objective function is characterized by a type parameter and information about the parameter is private. Specifically, each player only knows the range of the rivals' type parameters but not their true values.
Under such circumstances, 
rivals' type parameters are described by 
random variables and their actions are believed to depend 
on the type parameters.
  A general assumption is that the joint probability distribution of all type parameters
is publicly known.

This assumption allows each player to determine
the distribution of rivals' type parameters
conditional on its own type value and 
to take action to
optimize the conditional expected value of its objective function.
The Bayesian game is introduced by Harsanyi~\cite{harsanyi1967games} 
and has received much attention
in studies from modelling to equilibrium 
and has applications in strategic interactive decision-making problems with incomplete information, such as Bayesian network games \cite{blume2015linear,calvo2015communication}, 
rent-seeking games \cite{ewerhart2014unique,ewerhart2020unique,lockard2001efficient} and optimal public information disclosure problems \cite{colombo2014information,jia2023herding,morris2002social}.
For a comprehensive overview, 
see 
\cite{guo2021existence,meirowitz2003existence,ui2016bayesian} and
monograph~\cite{zamir2020bayesian}.
Classic Bayesian game models
focus on the case where the feasible set of each player
is independent of its type
and other players' actions. 
In practice, however, a player's feasible action
often depends on its rivals' actions
because it is constrained by
common resources and/or network capacities when they 
compete to provide homogeneous goods and services.
This kind of dependence is well-investigated in Nash games where 
the resulting games/equilibrium models are 
known as generalized Nash games/equilibrium models \cite{facchinei2010generalized}.
Likewise, a player's feasible action
might also depend 
on its own type. For example,
if we use the type parameter to denote a firm's marginal cost,
then its feasible action on pricing is inevitably 
subject to the marginal cost as discussed in \cite{liu2024Bayesian}. 
 Ewerhart and Quartieri \cite{ewerhart2020unique} 
 were the first to consider a Bayesian game for 
 imperfectly discriminating contests with finite type spaces
 where each player's feasible 
 action
 is subject to a type-dependent budget constraint. 
 They derived general conditions under which 
 a pure strategy Bayesian Nash equilibrium exists and is unique.

 In this paper, we follow the stream of research 
 to consider 
 a generalized Bayesian game where
 each player's feasible set depends on rivals' actions/strategies and its own type. 
Unlike \cite{ewerhart2020unique}, we consider
 type parameters that take continuous values and 
 investigate the existence and uniqueness 
 of a continuous pure strategy generalized Bayesian Nash equilibrium (GBNE).
 Following \cite{guo2021existence,meirowitz2003existence},
 we use Schauder's fixed point theorem to show 
  {\color{black}its} existence. A key step is to demonstrate the equicontinuity
of each player's optimal response function defined over the space of its type parameter. Unfortunately, 
 the dependence of  a player's feasible set on rivals' actions/strategies and  {\color{black} on} its own type makes this much more  {\color{black}complex}.
 Liu~et al. \cite{liu2024Bayesian} 
 considered 
a Bayesian Bertrand game with continuous type spaces
where each player's  {\color{black}feasible set of the} actions on pricing 
depends on its type (marginal cost) 
 {\color{black}and this class of dependence is modeled as the type-dependent interval for the feasible actions}. Under some moderate conditions, they derive the existence and uniqueness of a Bayesian equilibrium.
 {\color{black} Nevertheless, their model does not incorporate the dependence of each player's feasible set on its rivals' actions/strategies.}

The action-dependence of feasible sets  {\color{black}has been} 
explored  {\color{black}well} in the context of 
 {\color{black}the} generalized Nash game.
This class of games, which can be seen as the complete information counterpart of generalized Bayesian games, has been studied  {\color{black}extensively} and applied in various fields, such as economics \cite{arrow1954existence,tao2018duopoly},
computer science \cite{pang2008distributed}, 
 electric power distribution \cite{hobbs2007nash},
and environment management \cite{breton2006game}.
The existence and
uniqueness of
the generalized Nash equilibrium (GNE)
and numerical algorithms for computing such an equilibrium 
have been investigated using
fixed-point characterizations \cite{arrow1954existence,debreu1952social} and quasi-variational inequality (QVI) reformulations \cite{facchinei201012nash,harker1991generalized}.
We refer interested readers to the survey paper by Facchinei and Kanzow \cite{facchinei2010generalized} and Fischer et al. \cite{fischer2014generalized} for a comprehensive overview.
In recent years, there has been a growing interest in  stochastic GNE, see \cite{demiguel2009stochastic,lei2022stochastic,ravat2017existence,xu2013stochastic}, and GNE with infinite-dimensional action spaces, particularly in the context of optimal control where the players' strategy spaces can be Banach spaces \cite{kanzow2019multiplier} or Hilbert spaces \cite{borgens2021admm}. 
However, as mentioned above, since the best response strategy for each player in a (generalized) Bayesian game is a response 
function of its own type, many existing results for the (generalized) Nash games 
 {\color{black}can no longer} be
applied directly to the (generalized) Bayesian games. 


An important component of  {\color{black} this research} is to develop efficient numerical methods for computing 
an approximate GBNE. 
This is because generalized Bayesian games are infinite-dimensional,
 {\color{black} which makes it challenging to}
 {\color{black}determine}
an exact equilibrium, especially for the complex structures of type-and-action dependent feasible sets in the GBNE model. 
Athey \cite{athey2001single} investigated a finite-action step-like approximation for the non-decreasing pure strategy Bayesian Nash equilibrium.
Liu et al.~\cite{liu2024Bayesian} also developed a step-like approximation 
 {\color{black}for} the Bayesian Nash Bertrand equilibrium with type-dependent feasible sets and 
 {\color{black}applied}
a decomposition algorithm 
 {\color{black}to solve}
a
variational inequality (VI) problem reformulated from the Bayesian Nash equilibrium (BNE) model.
Another approach is to approximate the continuous Bayesian Nash equilibrium with polynomial functions. This method is widely applied to continuous linear programs \cite{bampou2012polynomial} and has been introduced to BNE models in \cite{guo2021existence}.
Here, we investigate the polynomial approximation scheme for the GBNE model and analyze the 
existence and convergence of polynomial GBNEs,
where 
additional effort is needed to deal with the
type-and-action dependent feasible set.

The main contributions of the paper can be summarized as follows.
}
\begin{enumerate}
    \item[(a)] \textbf{Novelty of the model}.
    We introduce a formal definition of the generalized Bayesian Nash equilibrium  with continuous type and action spaces, where the interdependence of players' feasible sets of actions is characterized by a 
    system of inequalities parameterized by their own types and rivals' response functions.
The new GBNE model 
subsumes the Bayesian Nash equilibrium  model with  type-dependent feasible sets (of actions). 
The latter 
 {\color{black}was}
recently proposed 
 {\color{black}
in \cite{ewerhart2020unique} 
in
the setting of a finite type space and investigated}
in \cite{liu2024Bayesian}
to describe the price competition
in discrete choice models,
 {\color{black} where the feasible set of each player’s action is a type-dependent interval.
This paper takes a step forward to 
generalize the existing 
BNE models by 
(a) including the interactions between players not only at the level of the objective functions but also at the level of the feasible sets by introducing 
type-and-action dependent feasible sets, and (b) 
proposing a unified framework where the feasible set is characterized by a system of inequalities, thus 
putting 
incorporating the Bayesian equilibrium (with/without the type-dependent feasible sets) and the generalized Nash equilibrium within a single framework.}
To illustrate why the GBNE model is useful in practice, we list three examples in rent seeking games, price competition games and 
Cournot 
games,  where  {\color{black}the} action space of each decision maker is subject to some type-dependent and/or action-dependent constraints. 

    
    \item[(b)] \textbf{Existence and uniqueness of continuous GBNE.}
We derive sufficient conditions 
 {\color{black}to establish} 
the existence of a continuous GBNE.
One of the main challenges to be tackled is the equicontinuity of the optimal response functions of players when their feasible sets of actions are interdependent. 
For this purpose, we consider the case 
 {\color{black}where}
each player's feasible set of actions is represented by a  system of inequalities, and each player's optimal decision making problem is a generic 
parametric program with inequality constraints.
Under 
some moderate conditions, we demonstrate the uniform H\"older continuity of  {\color{black}the set-valued mapping of feasible solutions and the set-valued mapping of the optimal solutions }
to the generic 
parametric program.
With this, we move on to
    use Schauder's fixed point theorem to
    show the existence of a GBNE under some additional conditions.  
     {\color{black}Given} a strictly monotone condition, we show  {\color{black}the} uniqueness of a continuous GBNE 
when each player's action space 
 {\color{black}only depends}
on its own type. 
In the case 
 {\color{black}where}
each player's action space is also dependent on rivals' actions,
we give a simple example to show that the uniqueness of GBNE is not guaranteed under standard
monotone conditions.
    

    \item[(c)] \textbf{Approximating the continuous GBNE 
    with SGNEs.}
    We use the well-known polynomial decision rule to 
    approximate each player's optimal response functions and subsequently 
    approximate the
    GBNE with 
    finite dimensional stochastic GNEs.
  This kind of 
  approximation method has been proposed by Guo et al.~\cite{guo2021existence}
for classic BNE models. 
When applied to the GBNE model,
new challenges
arise in
the type/action-dependent feasible sets 
where 
feasibility may be affected by the approximation.
We derive sufficient conditions to ensure the non-emptiness 
of the approximated feasible sets and 
consequently demonstrate the
 convergence of 
 polynomial GNEs to the GBNE as the degree of
  {\color{black}the}
 polynomials 
  {\color{black}increases} to infinity.
 Some preliminary numerical tests show that the approximation works well. 
 
\end{enumerate}

The rest of the paper is organized as follows. In section \ref{sec-model}, we give a formal definition of GBNE and showcase potential applications.
In section \ref{sec:existence}, we 
discuss the existence and uniqueness of a continuous 
GBNE.
In section \ref{sec:polynomial}, we propose
a polynomial approximation scheme 
 {\color{black}to compute}
an approximate GBNE and 
analyze the convergence of approximate GBNEs.
Finally, in section \ref{sec:numerical}, we present preliminary numerical test results to evaluate the proposed 
GBNE model and the polynomial approximation scheme.

\section{The model}
\label{sec-model}
We consider a generalized Bayesian game with $n$ players. 
For 
simplicity of notation, let $N :=\{1,\dots,n\}$.
For $i\in N$, we denote 
the type parameter 
 {\color{black}by $\theta_i \in \Theta_i \subset \R^{d_i}$}
with
$\Theta_i$ being a convex set,
and  $a_i \in \mathbb{A}_i \subset \R^{n_i}$ 
the action of player $i$.
A profile of types is a vector $\theta :=(\theta_1,\dots,\theta_n)\in 
\Theta:=\Theta_1\times \cdots\times \Theta_n$
and a profile of actions is a vector $a :=(a_1,\dots,a_n)\in  \mathbb{A}:=\mathbb{A}_1\times \cdots\times \mathbb{A}_n$.
By convention, 
we denote by $a_{-i}:=(a_1,\dots,a_{i-1},a_{i+1},\dots,a_n) \in \mathbb{A}_{-i}$ and $\theta_{-i}:=(\theta_1,\dots,\theta_{i-1},\theta_{i+1},\dots,\theta_n)\in \Theta_{-i}$ the vectors of actions
and types of all players except player $i$,  {\color{black} respectively}.
For $i\in N$, 
player $i$'s utility is represented by a continuous function 
$u_i(a_i,a_{-i},\theta)$ which maps from 
$\mathbb{A}_i\times \mathbb{A}_{-i} \times \Theta$ to $\R$.

In  Bayesian games, each player's type 
is 
private information, which means that
player $i$ only knows its own type $\theta_i$ but not its rivals' types $\theta_{-i}$.
To address the
incomplete information about the rivals' 
types, each player treats its rival $j$'s type parameter $\theta_j,j\neq i$ as a random variable, which means that $\theta_{-i}$ is a random vector from player $i$'s perspective. 
A standard assumption is that
the joint probability distribution of $\theta$, denoted by $\eta(\theta)$,
is public information.
We denoted the marginal distribution of $\theta_i$ for $i\in N$ by $\eta_i(\theta_i)$.
Conditional  on 
$\theta_i$,
player $i$'s posterior belief about 
the probability distribution 
$\theta_{-i}$ 
can be derived  {\color{black}as}
$\eta_i(\theta_{-i}|\theta_i) = \eta(\theta)/\eta_i(\theta_i)$,
which describes player $i$'s 
belief  {\color{black}about} its rivals' types.
In Bayesian games with
pure strategy, 
the optimal response function (strategy) 
of 
player $i$ is represented by 
$f_i$ which maps  from the type space $\Theta_i$ to the action space $\mathbb{A}_i$.
For $i\in N$, we denote by $\mathcal{F}_i$ the set of measurable functions $f_i: \Theta_i \rightarrow \mathbb{A}_i$ endowed with the infinity norm
$\|f_i\|_\infty = \max_{\theta_i\in \Theta_i}\|f_i(\theta_i)\|$. By convention, we use
$\|\cdot\|$ to denote the Euclidean norm in finite dimensional spaces, and 
$\mathcal{F} :=\prod_{j\in N} \mathcal{F}_j,
\mathcal{F}_{-i} := \prod_{j\in N\backslash \{i\}} \mathcal{F}_j$. 
Likewise, we denote 
the set of continuous functions $f_i:\Theta_i\rightarrow \mathcal{A}_i$ 
 {\color{black}by $\mathcal{C}_i$}, $\mathcal{C} :=\prod_{j\in N}\mathcal{C}_j, \mathcal{C} _{-i}:=\prod_{j\in N\backslash \{i\}}\mathcal{C}_j$.

\subsection{Generalized Bayesian Nash equilibrium}
In the framework of 
the \linebreak
Bayesian game, each player determines its optimal response strategy 
to 
maximize the expected value of its objective function w.r.t. the conditional distribution $\eta_i(\theta_{-i}|\theta_i)$. 
The best strategy for each player is a function $f_i(\cdot)$ on its type $\theta_i$.
In this paper, we consider the case that 
player $i$'s optimal action at its type $\theta_i$, written  $f_i^*(\theta_i)$, is obtained by solving the
following 
expected utility maximization problem
\begin{subequations}
\label{eqn:GBNE-player i}
  \bgeqn
        f_i^*(\theta_i) \in \arg \max\limits_{a_i\in \mathbb{A}_i} && \int_{\Theta_{-i}} u_i(a_i,f_{-i}(\theta_{-i}),\theta) d\eta_i(\theta_{-i}|\theta_i)\\
          \inmat{s.t.} &&
        \int_{\Theta_{-i}}g_i(a_i,f_{-i}(\theta_{-i}),\theta)d\eta_i(\theta_{-i}|\theta_i)\leq 0,    
    \label{eqn:GBNE-player i-b}
    \edeqn        
\end{subequations}
where $g_i(a_i,a_{-i},
\theta):\mathbb{A}_i\times\mathbb{A}_{-i}\times \Theta\rightarrow \R^{m_i}$ 
is a continuous function. 

The feasible set for each player $i$'s action at its type $\theta_i$
is defined 
 {\color{black}by}
the expected value 
of $g_i(a_i,f_{-i}(\theta_{-i}),\theta)$
which means that the effect from 
$f_{-i}(\theta_{-i})$ on the feasible set
is taken from an ``average/risk neutral'' point of view. It is possible to consider risk-averse or risk-taking constraints
by replacing $\mathbb{E}_{\eta_i(\theta_{-i}|\theta_i)}[g_i(a_i,f_{-i}(\theta_{-i}),\theta)]$ 
with
$R_i(g_i(a_i,f_{-i}(\theta_{-i}),\theta))$, where $R_i$ is a risk measure. In the case that $R_i$ is conditional value-at-risk (CVaR), we 
can still 
represent   $R_i(g_i(a_i,f_{-i}(\theta_{-i}),\theta))$ in the expectation form.
We will come back to this in Example \ref{exam:cournot_game}.
Throughout the paper, we make a blanket assumption that the integrals
 $$
 \int_{\Theta_{-i}} u_i(f_i(\theta_i), f_{-i}(\theta_{-i}),\theta)d\eta_i
(\theta_{-i}|\theta_i)\quad 
\inmat{and}\quad  \int_{\Theta_{-i}} g_i(f_i(\theta_i), f_{-i}(\theta_{-i}),\theta)d\eta_i (\theta_{-i}|\theta_i)
$$ 
are well-defined for each  fixed $\theta_i\in \Theta_i$, $f_i \in \mathcal{F}_i$ and $f_{-i}\in \mathcal{F}_{-i}$.
\begin{definition}[Generalized Bayesian Nash equilibrium (GBNE)] 
A pure strategy GBNE is an $n$-tuple $f^*:=(f_1^*,\dots,f_n^*)\in\mathcal{F}$ mapping from $\Theta_1\times\cdots\times\Theta_n$ to $\mathbb{A}_1\times\cdots\times\mathbb{A}_n$ such that for each $i\in N$ and fixed $\theta_i\in \Theta_i$,
    \bgeqn
    \begin{aligned}
        f_i^*(\theta_i) \in \arg \max\limits_{a_i\in \mathbb{A}_i} & \int_{\Theta_{-i}} u_i(a_i,f_{-i}^*(\theta_{-i}),\theta) d\eta_i(\theta_{-i}|\theta_i)\\
          \inmat{s.t.} &
        \int_{\Theta_{-i}}g_i(a_i,f_{-i}^*(\theta_{-i}),\theta)d\eta_i(\theta_{-i}|\theta_i)\leq 0 
        .    
    \end{aligned}
    \label{eqn:GBNE}
    \edeqn
\end{definition} 
 {\color{black} The main feature 
of this model which distinguishes 
it from 
existing BNE models is 
each player's feasible set of actions: 
here the feasible set depends 
not only on its only type parameter 
but also
on the rivals' response functions  
$f_{-i}$. Moreover, the feasible set is represented by a generic system of inequalities rather than 
by
an interval as in \cite{ewerhart2020unique,liu2024Bayesian}. 
The new feature 
enables the resulting 
GBNE
model to be 
 {\color{black}applied}
to a broader class of problems.
However, new challenges arise
when discussing the existence of an equilibrium.
We will 
return to this in Section 3.
}

Let
\begin{subequations}
\begin{eqnarray}
    \rho_i(a_i,f_{-i},\theta_i)& :=& \int_{\Theta_{-i}}u_i(a_i,f_{-i}(\theta_{-i}),\theta )d\eta_i(\theta_{-i}|\theta_i),\\
    h_i(a_i,f_{-i},\theta_i)& :=& \int_{\Theta_{-i}}g_i(a_i,f_{-i}(\theta_{-i}),\theta) d\eta_i(\theta_{-i}|\theta_i).
\end{eqnarray}
    \label{eqn:def_of_h}
\end{subequations}
Then we can write the 
GBNE model (\ref{eqn:GBNE}) succinctly as
\bgeqn
    \begin{aligned}
        f_i^*(\theta_i)\in\arg \max\limits_{a_i\in \mathbb{A}_i} &\; \rho_i(a_i,f_{-i}^*,\theta_i), \\
        \inmat{s.t.} &\; h_i(a_i,f_{-i}^*,\theta_i)\leq 0, 
    \end{aligned}
    \label{eqn:GNE} 
\edeqn
for each $i\in N$, $\theta_i\in \Theta_i$.
This is an infinite-dimensional 
equilibrium problem with the best response strategies  {\color{black}being} defined over the type space $\Theta_i$, rather than  {\color{black}by} a fixed action in the generalized Nash game in \cite{facchinei2010generalized}.
In the case that the type space $\Theta$ shrinks to a singleton, model (\ref{eqn:GNE}) reduces to a classic 
generalized Nash equilibrium 
model. 

For each fixed $f\in \mathcal{F}$, 
let
 \bgeqn
     \mathcal{A}_i(f_{-i},\theta_i) := \left\{ a_i\in \mathbb{A}_i: \int_{\Theta_{-i}}g_i(a_i,f_{-i}(\theta_{-i}),\theta)d\eta_i(\theta_{-i}|\theta_i)\leq 0 \right\}
    \label{def:feasible_set}
 \edeqn
and define
\bgeqn \label{eqn:mapping}
    \Psi (f)&:=& \bigg\{ (\tilde{f}_1,\dots,\tilde{f}_n)\in \mathcal{F}:
    \inmat{for} \; i\in N, 
    \\
    && \tilde{f}_i(\theta_i)\in \arg \max_{a_i\in \mathcal{A}_i(f_{-i},\theta_i)} \rho_i(a_i,f_{-i},\theta_i), \quad \forall \theta_i\in \Theta_i\bigg\}.\nonumber
\edeqn
Then $f^*$ is a GBNE iff 
$f^* \in \Psi (f^*)$.
In the case when the GBNE is continuous, problem (\ref{eqn:GBNE}) can be represented in terms of the expectation of the utility functions with respect to  distribution $\eta$, and subsequently be
reformulated as a single optimization problem. The next proposition states this.


\begin{proposition} [Equivalent formulations of the GBNE model] \label{thm:equivalence}
     Assume  for $i\in N$ and any fixed $(a_{-i},\theta)$, $u_i(a_i,a_{-i},\theta)$ is continuous in $a_i$, and $g_i(a_i,a_{-i},\theta)$  is convex in $a_i$.
     Then $f^*$ is a continuous
     GBNE of model (\ref{eqn:GBNE}) iff 
     \begin{itemize}
         \item[(i)] 
          for $i\in N$,
          \bgeqn\label{eqn:GBNE'2}
          \begin{aligned}
              f_i^*\in \arg  \max_{f_i\in \mathcal{C}_i} &\;\mathbb{E}_{\eta} [u_i(f_i(\theta_i),f_{-i}^*(\theta_{-i}),\theta)]\\
                \inmat{s.t.} &\; h_i(f_i(\theta_i),f_{-i}^*,\theta_i) \leq 0,\quad \forall \theta_i\in\Theta_i,
          \end{aligned}
            \edeqn
        \item[(ii)] or 
        equivalently
            \bgeqn\label{eqn:GBNE''2}
                f^*\in \arg & \max_{f\in \mathcal{C}} &\sum_{i=1}^n \mathbb{E}_{\eta} [u_i(f_i(\theta_i),f_{-i}^*(\theta_{-i}),\theta)]\\\nonumber
                & \inmat{s.t.} & h_i(f_i(\theta_i),f_{-i}^*,\theta_i)\leq 0,\quad \forall \theta_i\in\Theta_i, i\in N.
            \edeqn
     \end{itemize} 
\end{proposition}

\noindent
\textbf{Proof.} 
Part (ii) follows directly from Part (i).
We only prove the latter.

 {\color{black}\textbf{The ``only if'' part.}} 
Let $f^*\in\mathcal{C}$ be a GBNE.
Then for
any 
$f_i\in \mathcal{C}_i,i\in N$ satisfying the constraints in model (\ref{eqn:GBNE'2}),
\bgeq
\int_{ \Theta_{-i}} u_i(f_i^*(\theta_i), f_{-i}^*(\theta_{-i}),\theta)d\eta_i(\theta_{-i}|\theta_i)\geq \int_{ \Theta_{-i}} u_i(f_i(\theta_i), f_{-i}^*(\theta_{-i}),\theta)d\eta_i(\theta_{-i}|\theta_i), \forall \theta_i \in \Theta_i.
\edeq
By integrating with respect to $\theta_i$ on both sides of the inequality, we immediately 
obtain (\ref{eqn:GBNE'2}).

 {\color{black}\textbf{The ``if'' part.}} Let $f^* \in \mathcal{C}$ and $f^*$ satisfy (\ref{eqn:GBNE'2}). We show that $f^*$ is a GBNE of model (\ref{eqn:GBNE}). Assume, for the sake of a contradiction, that $f^*$ is not a GBNE. Then there exists some $i\in N$ and $\hat{f}_{i}\in \mathcal{C}_i$ satisfying $h_i(\hat{f}_{i}(\theta_i),f^*_{-i},\theta_i)\leq 0, \forall \theta_i\in \Theta_i$ such that for some $\hat{\theta}_{i} \in \Theta_{i}$, 
\bgeq
\int_{\Theta_{-i}} u_i(f_i^*(\hat{\theta}_i), f_{-i}^*(\theta_{-i}),\hat{\theta})d\eta_i(\theta_{-i}|\hat{\theta}_i)
<\int_{\Theta_{-i}} u_i(\hat{f}_i(\hat{\theta}_{i}), f_{-i}^*(\theta_{-i}),\hat{\theta})d\eta_i(\theta_{-i}|\hat{\theta}_{i}).
\edeq
Then, the inequality above implies that there exists a neighborhood of $\hat{\theta}_i$, written $\mathcal{B}_{\hat{\theta}_i}\subset \Theta_i$,  such that
\bgeqn
\label{eqn:equivalent_proof}
    & &\int_{\mathcal{B}_{\hat{\theta}_i}}\int_{\Theta_{-i}} u_i(f_i^*(\theta_i), f_{-i}^*(\theta_{-i}),\theta_i,\theta_{-i})d\eta_i(\theta_{-i}|\theta_i)d \eta_i(\theta_i)\\\nonumber
    &<& \int_{\mathcal{B}_{\hat{\theta}_i}}\int_{\Theta_{-i}} u_i(\hat{f}_i(\theta_{i}), f_{-i}^*(\theta_{-i}),\theta_{i},\theta_{-i})d\eta_i(\theta_{-i}|\theta_{i}) d \eta_i(\theta_i).    
\edeqn
Consequently, we can construct a continuous function $\tilde{f}_i(\theta_i)$ such that $\tilde{f}_i(\theta_i) =
f_i^*(\theta_i)$ for $\theta_i \notin \mathcal{B}_{\hat{\theta}_i}$
and $\tilde{f}_i(\theta_i) =
\hat{f}_i(\theta_i)$ otherwise. Note that in general 
the function values of $\hat{f}_i$ do not necessarily 
meet those of $f^*_i$ at the boundaries of $\mathcal{B}_{\hat{\theta}_i}$, in which case, we may revise the function values of $\tilde{f}_i$ near the boundaries within 
$\mathcal{B}_{\hat{\theta}_i}$ such that the inequality (\ref{eqn:equivalent_proof}) is preserved.
Next, we can write $\tilde{f}_i(\theta_i)$ as 
\bgeq
\tilde{f}_i(\theta_i):= \lambda(\theta_i)f_i^*(\theta_i) + (1-\lambda(\theta_i)) \hat{f}_i(\theta_i), \forall \theta_i\in \Theta_i,
\edeq
where $\lambda(\theta_i)=1$ for $\theta_i \notin \mathcal{B}_{\hat{\theta}_i}$ and $\lambda(\theta_i)=0$ otherwise. 
By the convexity of \linebreak $h(a_i,f_{-i}(\theta_{-i}),\theta_i)$ in $a_i$, 
$$
h(\tilde{f}_{i}(\theta_i),f_{-i}^*,\theta_i)\leq \lambda(\theta_i) h(f_{i}^*(\theta_i),f_{-i}^*,\theta_i) + (1-\lambda(\theta_i)) h(\hat{f}_{i}^*(\theta_i),f_{-i}^*,\theta_i) \leq 0, \forall \theta_i\in \Theta_i.
$$
On the other hand,
\bgeq
    \mathbb{E}_\eta[u_i(f_i^*(\theta_i),f_{-i}^*(\theta_{-i}),\theta)] < \mathbb{E}_\eta[u_i(\tilde{f}_i(\theta_i),f_{-i}^*(\theta_{-i}),\theta)],
\edeq 
which leads to a contradiction to (\ref{eqn:GBNE'2}) as desired.
\hfill $\Box$


In the case 
 {\color{black}where}
the objective functions are continuously differentiable,
we can use model (\ref{eqn:GBNE'2})
to characterize the GBNE model (\ref{eqn:GBNE}) as 
the following 
stochastic quasi-variational inequality (QVI) problem:
\bgeqn
    \mathbb{E}_\eta [\nabla_a u(f^*(\theta),\theta_i)^T(\tilde{f}(\theta)-f^*(\theta))]\leq 0, \quad\forall \tilde{f}\in \mathcal{C}(f^*),
    \label{eqn:qvi-GBNE}
\edeqn
where $u(f(\theta),\theta_i) := (u_1(f_1(\theta_1),f_{-1}(\theta_{-1}),\theta_1),\dots,u_n(f_n(\theta_n),f_{-n}(\theta_{-n}),\theta_n))^T$ and \linebreak
\bgeqn \label{def:cont_func_set}
\mathcal{C}(f):= \{\tilde{f}\in \mathcal{C}: h_i(\tilde{f}_i(\theta_i),f_{-i},\theta_i)\leq 0, \forall \theta_i\in \Theta_i, i\in N\}
\edeqn
is a feasible set of continuous response functions. 
The terminology ``quasi-variational'' is used to 
highlight
the dependence of the feasible set on $f^*$, see \cite{facchinei2010generalized,lei2022stochastic}.
When the feasible set is independent of
$f^*$, 
problem (\ref{eqn:qvi-GBNE}) reduces to 
an infinite-dimensional variational inequality problem
via which  Ui \cite{ui2016bayesian} 
explored the existence and uniqueness of the classic BNE.

\subsection{Applications}
\label{subsec:application}
Before moving on to 
theoretical discussion of the GBNE model,  
we 
give three examples 
to illustrate how the model may be applied in practice. 

\begin{example} [Rent-seeking games]
    Rent-seeking games \cite{lockard2001efficient} with budget constraints
    have been considered in  {\color{black}the} literature, see e.g.~\cite{che1998caps,ewerhart2020unique,leininger1991patent}. 
    The type-dependent budget constraints for rent-seeking games with finite type spaces have been investigated in \cite{ewerhart2020unique}. 
    Here we consider a slightly more general case where the type parameter takes {\em continuous values} to fit 
    our GBNE model.
    In this game, $n$ players choose levels of the costly effort in order to win a share of the contested rent, and each player 
     {\color{black}has} 
    private information about the data of the game.
    Suppose that player $i$'s cost is a linear function $c_i(a_i)=\theta_i a_i$ of its effort parameterized by the type variable $\theta_i$, and 
     {\color{black}the player $i$'s} budget constraint 
    $\mathcal{A}_i(\theta_i)=[0,a_i^\text{max}(\theta_i)]$ 
    is also parameterized by type $\theta_i$. The parameter $\theta$ is  {\color{black} assumed to be} drawn from a joint probability distribution $\eta$, which is absolutely continuous w.r.t. the Lebesgue measures over its support $\Theta$. Suppose that player $i$'s utility function is $u_i(a,\theta) = c_i(a) + \frac{a_i}{a_i+\sum_{j\neq i}a_j}$ and  $\frac{a_i}{a_i+\sum_{j\neq i}a_j}$ is set to $1/n$ if $a_i+\sum_{j\neq i}a_j = 0$. 
    Then the rent-seeking game with 
     {\color{black}a}
    type-dependent budget set is to find an $n$-tuple $(f_1,\dots,f_n)$ such that for all $i\in N$, 
    \bgeqn
        f_i(\theta_i) \in \arg \max_{a_i\in \mathcal{A}_i(\theta_i)} -a_i\theta_i + \int_{\Theta_{-i}} \frac{a_i}{a_i+\sum_{j\neq i}f_j(\theta_j)} d\eta_i(\theta_{-i}|\theta_i), \forall \theta_i\in \Theta_i.
    \edeqn
    This is a special GBNE with the feasible set of actions depending  {\color{black}only} on each player's type.
\end{example}

\begin{example}[The price competition in discrete choice models \cite{liu2024Bayesian}]
 Consider an 
  {\color{black}oligopolistic}
market with $n$ differentiated 
substitutable  products where each product 
is offered by a distinct firm.
Each customer chooses to buy at most one product and their utility of purchasing product $j$ is determined by the following equation:
\bgeq 
U_j = \beta^{T} x_j - \alpha p_j + \epsilon_j, \quad \inmat{for} \;  j=1, \ldots, n,
\edeq
where $\beta$ is a vector representing  {\color{black}the} customer's preference or taste for the observed product characteristic $x_j$, $\alpha$ is a positive scalar representing customer's sensitivity to price $p_j$, and $\epsilon_j$ denotes the utility of unobserved characteristics or idiosyncratic product-specific random demand shock, 
with $\epsilon_1,\dots,\epsilon_n$ being independent and identically 
distributed and following the Gumbel distribution. 
Liu et al.~\cite{liu2024Bayesian} consider the 
price competition model
 {\color{black} to find an $n$-tuple $(f_1,\dots,f_n)$ such that for all $i\in N$, }
    \bgeqn
    \label{eqn:price_competition}
       \qquad  f_i(\theta_i) \in \arg \max\limits_{p_i\in \mathcal{A}_i(\theta_i)} \int_{\Theta_{-i}} (p_i-\theta_i)\frac{e^{\beta^Tx_i-\alpha p_i}}{1+\sum_{k=1}^n e^{\beta^Tx_k-\alpha p_k}} d\eta_i(\theta_{-i}|\theta_i), \forall \theta_i\in\Theta_i,
    \edeqn
where $\theta_i$ represents firm $i$'s marginal cost, 
and $\frac{e^{\beta^Tx_i-\alpha p_i}}{1+\sum_{k=1}^n e^{\beta^Tx_k-\alpha p_k}}$ is the probability that product $i$ is chosen by customers which is also known as the market share. In this model, the marginal cost of a firm is private information and hence treated 
as a random variable from rivals' perspectives. The joint probability distribution of the marginal costs is assumed 
to be 
 {\color{black}publicly known.}
A unique feature of the model is that
the feasible set $\mathcal{A}_i(\theta_i)$ 
is dependent on $\theta_i$ because the sale price must exceed the marginal cost under normal circumstances (without predatory pricing behavior). The price competition model (\ref{eqn:price_competition}) is a GBNE model.

\end{example}

\begin{example}[Cournot games]
\label{exam:cournot_game}
    Cournot games
    have been widely considered in the context of (stochastic) generalized Nash games, see e.g.~\cite{hobbs2007nash,ravat2017existence,tao2018duopoly}. 
    Here we investigate Cournot games with incomplete information which 
 fits the framework of the GBNE model.
    
    Consider a Cournot game with $n$ firms,
    each firm $i$ determines the quantity $a_i$ to produce in order to maximize 
    its expected payoff.
    The production 
    has a 
    capacity constraint $a_i^{\text{max}}$.
The cost function, written 
$c_i(a_i, \theta_i)$,
is 
continuous in $a_i$ for any $\theta_i\in \Theta_i$ and
characterized by 
a type parameter 
$\theta_i$ where
the joint distribution of $\theta:=(\theta_1,\dots,\theta_n)$, written  
$\eta$, is publicly known. 
The payoff function
is 
$u_i(a, \theta) := - c_i(a_i, \theta_i) + p (\sum_{j=1}^n a_j)a_i$, where $p(\cdot)$ is 
a continuous
inverse demand function.
 {\color{black}Unlike} the Cournot game in \cite{ui2016bayesian}, here
  $a_i$ must satisfy the following constraints
    \bgeqn 
    \label{eqn:coupled_cons}
    \tilde{g}_i(a_i,f_{-i}(\theta_{-i}) ,\theta)
    \leq 0,\quad \forall \theta_{-i}\in \Theta_{-i},
    \edeqn
    where $\tilde{g}_i(a_i,a_{-i},\theta):\mathbb{A}_i\times\mathbb{A}_{-i}\times \Theta \rightarrow \R$ is a continuous function.
    In practice,  {\color{black}these} constraints may be interpreted
    as the dependence of the availability of 
    resources 
     {\color{black}to}
    firm $i$ 
    on its rivals' actions; see \cite{facchinei2010generalized}.  
    Incorporating 
    constraints (\ref{eqn:coupled_cons}) 
    directly into a Bayesian Nash equilibrium model may 
    lead to a semi-infinite 
    program which is difficult to handle.
    One way to tackle the difficulty is to 
    set 
    $g_i$ in GBNE model (\ref{eqn:GBNE}) as $g_i(a_i,f_{-i}(\theta_{-i}),\theta) := \tilde{g}_i(a_i,f_{-i}(\theta_{-i}),\theta)$ and 
    replace the semi-infinite constraints with a single 
expected value constraint:
    \bgeqn
    \label{eqn:expected_cons}
        \mathbb{E}_{\eta(\theta_{-i}|\theta_i)} [g_i(a_i,f_{-i}(\theta_{-i}),\theta)] := \mathbb{E}_{\eta(\theta_{-i}|\theta_i)} [\tilde{g}_i(a_i,f_{-i}(\theta_{-i}),\theta)]\leq 0.
    \edeqn
    However, (\ref{eqn:expected_cons})  has a serious drawback because it significantly relaxes constraints
   (\ref{eqn:coupled_cons}). 
  
    Alternatively, we may 
    replace the semi-infinite
    constraints (\ref{eqn:coupled_cons}) by a chance constraint  
    $
        \mathbb{P}_{\theta_{-i}|\theta_i}\left(\tilde{g}(a_i,f_{-i}(\theta_{-i}),\theta)\leq 0\right)  \geq 1-\beta 
    $,
    where $\beta\in (0,1)$, see  e.g.~\cite{campi2018wait}. 
    Under some mild conditions, the chance constraint may be 
    approximated by 
    the conditional Value-at-Risk (CVaR) constraint   {\color{black} w.r.t. $\theta_{-i}$ conditional on $\theta_i$ (see e.g.~\cite{nemirovski2007convex}):
    \begin{eqnarray} \label{eqn:cvar}
        \text{CVaR}_{1-\beta} (\tilde{g}(a_i,f_{-i}(\theta_{-i}),\theta)|\theta_i)\leq 0,
    \end{eqnarray}
    where 
    $
        \text{CVaR}_{1-\beta}(X):= \frac{1}{\beta}\int_{X\geq\text{VaR}_{1-\beta}(X)} x d F_X(x)
    $,
    $F_X(x)$ denotes the cumulative distribution function of a random variable $X$, and  $\text{VaR}_{1-\beta}(X) := \inf \{\alpha\in \R:\mathbb{P}(X\leq \alpha) \geq 1-\beta \ \} $ is the value-at-risk (VaR).
    To ease the discussion, we make a blanket assumption that for any $\theta_i\in \Theta_i$ and $f_{-i}\in \mathcal{F}_{-i}$, there exists a feasible action $a_i$ 
    satisfying the constraint (\ref{eqn:cvar}).
    By the Rockafellar/Uryasev reformulation (see \cite{rockafellar2000optimization}), 
    \begin{eqnarray*}
         \text{CVaR}_{1-\beta}(X) = \min_{t\in \R} t+\frac{1}{\beta} \int_{\text{supp}(X)} (X-t)_+ dF_X(x),
    \end{eqnarray*}
    where $(x)_+ = \max \{x,0\}$ and $\text{supp}(X)$ denotes the support of $X$.
    Using the formulation, we can recast 
    the GBNE model 
    with the CVaR constraints (\ref{eqn:cvar}) as: 
    finding an $n$-tuple $(f_1,\dots,f_n)$ such that for all $i\in N$, 
    \begin{eqnarray}
    \begin{aligned}
        \quad f_i(\theta_i) \in \arg  \max\limits_{a_i \in [0,a_i^{\text{max}}]} & - c_i(a_i, \theta_i) + a_i\int_{\Theta_{-i}} p(a_i+\sum_{j\neq i} f_j(\theta_j))  d\eta_i(\theta_{-i}|\theta_i)\\
        \text{s.t.} & \min_{t_i\in \R} t_i + \frac{1}{\beta}\int_{\Theta_{-i}} (\tilde{g}_i(a_i,f_{-i}(\theta_{-i}),\theta) - t_i)_+ d\eta_i(\theta_{-i}|\theta_i) \leq 0, \label{eqn:cournot_min_constraint}
    \end{aligned}
    \end{eqnarray}
    for 
    $\theta_i\in \Theta_i$.
     The problem is well defined in the sense that 
     the optimum is finite valued and is attainable.
    By making $t_i$ 
    a variable, we can further reformulate the problem as:
    finding
    an $n$-tuple $(f_1,\dots,f_n)$ such that for 
    $i\in N$ and each fixed $\theta_i\in \Theta_i$
    \bgeqn
    \begin{aligned} 
    (f_i(\theta_i),t_i(\theta_i))
        \in \arg  \max_{a_i \in [0,a_i^{\text{max}}],t_i\in \R} & - c_i(a_i, \theta_i) + a_i\int_{\Theta_{-i}} p(a_i+\sum_{j\neq i} f_j(\theta_j))  d\eta_i(\theta_{-i}|\theta_i)\\
        \inmat{s.t.}\; & t_i + \frac{1}{\beta}\int_{\Theta_{-i}} (\tilde{g}_i(a_i,f_{-i}(\theta_{-i}),\theta) - t_i)_+ d\eta_i(\theta_{-i}|\theta_i) \leq 0
    \end{aligned}
    \label{eq:GBNE-ex-CVaR}
    \edeqn
      for some $t(\theta_i)\in \R$.
Equivalence between 
(\ref{eqn:cvar})
and (\ref{eq:GBNE-ex-CVaR})
can be established by showing 
that the optimal solution of the former can be used to construct a feasible solution of the latter and vice versa,  consequently the optimal values are equal given that the objective functions are identical.
     Furthermore, by setting $g_i(a_i,f_{-i}(\theta_{-i}),\theta) := t_i + \frac{1}{\beta}(\tilde{g}_i(a_i,f_{-i}(\theta_{-i}),\theta) -t_i)_+,$
    we can fit 
    (\ref{eq:GBNE-ex-CVaR})
    into the GBNE framework 
    (\ref{eqn:GBNE}). 
    }
\end{example}

\section{Existence and uniqueness of GBNE}\label{sec:existence}
In this section, we 
derive sufficient conditions
under which problem (\ref{eqn:GBNE}) has  
a  GBNE.
A key step towards 
establishing the existence of 
a GBNE is to derive the equicontinuity of the 
optimal response function of each player. 
For this purpose, we 
consider a generic parametric program corresponding to a player's optimal decision-making problem and  
investigate 
the uniform continuity of 
the optimal solution set mapping of 
this program with respect to  the parameters.


\subsection{
 Uniform continuity of the
 optimal solution set mapping of a parametric program}
 \label{subsec:continuity_of_para_program}
Consider a parametric minimization problem
\bgeqn
    \label{eqn:para_min_problem}
    \min_{x\in \mathbb{X}} \quad &&  \phi(x,y,t), \nonumber\\
 \text{s.t.} \quad &&
 \psi_i(x,y,t)\leq 0, \; \inmat{for} \;  i = 1,\dots, m,
\edeqn
where  $\mathbb{X}\subset\R^n$, $\mathbb{Y}$ is a Banach space equipped with the norm $\|\cdot\|$; 
 {\color{black}
$Y\subset \mathbb{Y}$ is a nonempty subset and  
$T\subset \R^l$ is a nonempty and convex subset; 
$(y,t)\in  Y\times T$ are parameters;
}
$\phi: \mathbb{X}\times Y \times T \rightarrow \R$ and 
$\psi_i:\mathbb{X}\times Y \times T \rightarrow \R$, $i = 1,\dots, m$
are continuous functions. 
To ease the exposition, let $\psi(x,y,t):= (\psi_1(x,y,t),\dots,\psi_m(x,y,t))^T$ and  $X(y,t)\subset \mathbb{X}$ denote the feasible set. Let $\phi^*(y,t)$ and  $X^*(y,t)$ 
denote the optimal value and the set of optimal solutions respectively. 
Our aim here is to derive the equicontinuity of 
$X^*(y,t)$ in $t$ {\color{black}:} for any $\epsilon>0$, there exist $\delta>0$ such that
\bgeqn
\label{eqn:equicontinuity}
\sup_{y\in Y} 
\mathbb{H}(X^*(y,t'),X^*(y,t''))
< \epsilon, \forall t',t''\in T, \; \inmat{with} \; \|t'-t''\| \leq \delta.
\edeqn
  Here and  {\color{black}below}
  $\mathbb{H}$ denotes the Hausdorff distance 
        between two non-empty  compact sets.
To this end, we  {\color{black} first prove}
the feasible set-valued mapping $X(y,\cdot)$ to be 
H\"older continuous in $t$
 {\color{black} in the following lemma, and then followed  
by 
the equicontinuity of $X^*(y,t)$ in $t$ 
in the subsequent theorem.}


\begin{lemma}[Uniform continuity of the feasible set of (\ref{eqn:para_min_problem})] \label{lemma:hoffman}
    Consider the feasible set of problem (\ref{eqn:para_min_problem}). 
    Assume: (a) $\psi_i(x,y,t)$, $i=1,\dots,m$ is convex in $x$ for all fixed $(y,t)$; 
    (b) 
    Slater's condition holds
    uniformly
    for all $(y,t)\in Y\times T$ 
    with a positive constant $\alpha$ {\color{black}:}
    for any $(y,t)$, 
    there exists a feasible point, denoted by $\bar{x}(y,t)
    \in X(y,t)$, 
    satisfying
    \bgeqn
        \psi_i(\bar{x}(y,t),y,t)\leq -\alpha,\ \inmat{for} \;  i = 1,\dots, m; 
        \label{eqn:slater}
    \edeqn
    (c) there exists a positive constant $A$ such that
        $\inmat{diam}(X(y,t))\leq A$ for all $(y,t)$, 
        where $\inmat{diam}(X) := \sup_{x,x'\in X} \|x-x'\|$;
    (d) $\psi_i(x,y,t)$ is H{\"o}lder continuous in $t\in T$ 
        uniformly for $(x,y)\in \mathbb{X}\times Y$ {\color{black}:}
        there exist positive constants $\kappa_\psi$ and $q_\psi\in (0,1]$ such that
        \bgeq
            \left\| \psi(x,y,t')-\psi(x,y,t'') \right\| \leq \kappa_\psi \| t' - t''\| ^{q_\psi}, \forall t', t''\in T, 
        \edeq
       uniformly for $(x,y)\in \mathbb{X}\times Y$.
    Then the following assertions hold.
    \begin{itemize}
        \item[(i)] 
        For any fixed $(y,t)\in Y\times T$ and
        $x \in \R^n$,
     let $\gamma:= \|(\psi(x,y,t))_+\|$ and $\hat{x}(y,t):= (1-\frac{\gamma}{\gamma+\alpha})x + \frac{\gamma}{\gamma+\alpha} \bar{x}(y,t)$. Then  $\hat{x}(y,t)\in X(y,t)$,
        \bgeqn \label{eq:Hoff-Lem-1}
            \|x-\hat{x}(y,t)\|
            &\leq&
                        \frac{\|x-\bar{x}(y,t)\|}{\alpha}\|(\psi(x,y,t))_+\|
        \edeqn
        and
        \bgeqn
            \dd(x,X(y,t))&\leq&
            \frac{\inmat{diam}(X(y,t))}{\alpha}\|(\psi(x,y,t))_+\|, 
            \label{eq:Hoff-Lem-2}
        \edeqn 
        where $(x)_+ := \max(0,x)$ with the maximum being taken componentwise, and $\dd(x,X):= \inf_{\tilde{x}\in X}\|x-\tilde{x}\|$.
        \item[(ii)] The feasible  {\color{black}set-valued} mapping $X(y,t)$ is uniformly 
        continuous in $t\in T$: for all $y\in Y$, 
        \bgeqn
        \label{eqn:stability_of_feasible_set}
            \mathbb{H}(X(y,t'),X(y,t'')) 
            \leq \frac{A\kappa_\psi} {\alpha} \| t' - t''\|^{q_\psi}, \forall t',t''\in T.
        \edeqn

        \item[(iii)] For any fixed $y\in Y$, there exists a continuous mapping $\upsilon:T\rightarrow \mathbb{X}$ such that for $i=1,\dots,m$,
\bgeqn
\label{eqn:hoffman_corollary}
    \psi_i(\upsilon(t),y,t)\leq -\frac{\alpha}{2},\; \forall t\in T.
\edeqn
      

    \end{itemize}
\end{lemma}

\noindent
\textbf{Proof.}  Part (i).
Let $(y,t)\in Y\times T$ and $x\in \R^n$ be fixed and $\gamma:= \|(\psi(x,y,t))_+\|$. Then $\psi(x,y,t)\leq \gamma \mathbf{e}$, where $\mathbf{e}$ is an $m$-dimensional vector with unit components.
By the convexity of $\psi(x,y,t)$ in $x$ and the definition of $\hat{x}(y,t)= (1-\frac{\gamma}{\gamma+\alpha})x + \frac{\gamma}{\gamma+\alpha} \bar{x}(y,t)$, we have
\bgeq
    \psi(\hat{x}(y,t),y,t) & \leq & \left(1-\frac{\gamma}{\gamma+\alpha}\right) \psi(x,y,t) + \frac{\gamma}{\gamma+\alpha} \psi(\bar{x}(y,t),y,t)\\
    & \leq & \left(1-\frac{\gamma}{\gamma+\alpha}\right) \gamma\mathbf{e} + \frac{\gamma}{\gamma+\alpha} (-\alpha\mathbf{e}) = 0,
\edeq
where the inequality is taken componentwise.
This shows $\hat{x}(y,t)\in X(y,t)$. 
By the definition of $\hat{x}(y,t)$,
\bgeq 
     \|x-\hat{x}(y,t)\|
     = \frac{\gamma}{\gamma+\alpha} \| x - \bar{x}(y,t) \| 
     \leq \frac{\|(\psi(x,y,t))_+\|}{\alpha}\| x - \bar{x}(y,t) \|,
\edeq
which 
yields (\ref{eq:Hoff-Lem-1}).
Observe that
$
    \alpha  \|x-\hat{x}(y,t)\| = \gamma (\| x - \bar{x}(y,t) \| - \|x-\hat{x}(y,t)\|).
$
Thus
\bgeq
    \dd(x,X(y,t)) &\leq& \|x-\hat{x}(y,t)\| = \frac{\gamma}{\alpha} (\| x - \bar{x}(y,t) \| - \|x-\hat{x}(y,t)\|)\leq \frac{\gamma}{\alpha} \| \bar{x}(y,t) - \hat{x}(y,t) \| \\
    &\leq& \frac{\|(\psi(x,y,t))_+\|}{\alpha} \inmat{diam}(X(y,t)),
\edeq
which is (\ref{eq:Hoff-Lem-2}).

Part (ii). For any $y\in Y$, let $t',t''\in T$, and $x'\in X(y,t')$. By Part (i),  {\color{black} condition (c)} and the H{\"o}lder continuity in condition (d),
\bgeq
    \dd(x',X(y,t'')) & \leq & \frac{\inmat{diam}(X(y,t''))}{\alpha}\|(\psi(x',y,t''))_+\|\\
    & \leq & \frac{A}{\alpha} \left\|(\psi(x',y,t''))_+ - (\psi(x',y,t'))_+\right\|\\
    & \leq & \frac{A}{\alpha} \left\| \psi(x',y,t'')-\psi(x',y,t') \right\|\\
    & \leq & \frac{A\kappa_\psi}{\alpha}  \| t' - t''\| ^{q_\psi},
\edeq
 {\color{black}where, by the definition of $X(y,t')$, $x'\in X(y,t')$ implies that $\psi(x',y,t')\leq 0$ and thus $(\psi(x',y,t'))_+= 0$.} 
This 
 {\color{black} implies} that for all $y\in Y$,
\bgeqn
\label{eqn:D_1_2}
    \mathbb{D}(X(y,t'),X(y,t''))
    \leq \frac{A\kappa_\psi}{\alpha}\| t' - t''\| ^{q_\psi}, \forall t',t''\in T,
\edeqn
where $\mathbb{D}(S_1,S_2) :=\sup_{s\in S_1}\dd(s,S_2)$ for two sets $S_1, S_2$.
By swapping the positions between $X(y,t')$ and $X(y,t'')$, we obtain that for all $y\in Y$,
\bgeqn
\label{eqn:D_2_1}
    \mathbb{D}(X(y,t''),X(y,t'))
    \leq \frac{A\kappa_\psi}{\alpha} \| t'' - t'\| ^{q_\psi}, \forall t',t''\in T.
\edeqn
A combination of (\ref{eqn:D_1_2}) and (\ref{eqn:D_2_1}) 
gives rise to 
(\ref{eqn:stability_of_feasible_set}).

Part (iii). 
Consider the following parametric 
system of 
inequalities
\bgeqn
\label{eqn:hoffman_corollary_slater}
\psi_i(x,y,t) + \frac{\alpha}{2} \leq 0, \; \inmat{for} \;  i = 1,\dots, m.
\edeqn
Observe that (\ref{eqn:hoffman_corollary_slater}) satisfies 
all 
the conditions of Lemma \ref{lemma:hoffman} except that the constant $\alpha$ in Slater's condition is replaced with
$\frac{\alpha}{2}$ in condition (c). 
By Lemma \ref{lemma:hoffman} (ii), 
the set of solutions to  the system 
(\ref{eqn:hoffman_corollary_slater}), written $\bar{X}(y,t)$,
is  convex set-valued 
and
H\"older continuous in $t$ over $T$.
The conclusion (\ref{eqn:hoffman_corollary}) follows from 
the continuous selection theorem (see e.g. Theorem 9.1.2 in \cite{aubin2009set}) 
which ensures
that
$\upsilon(\cdot)$ 
 {\color{black}is} continuous.
\hfill $\Box$

The next theorem quantifies the stability of
the set of optimal solutions $X^*(y,t)$ when parameter $t$ varies and then derives equicontinuity of $X^*(y,t)$ in $t$.

\begin{theorem}
[Uniform continuity of the optimal solution
set of (\ref{eqn:para_min_problem})]
\label{thm:stability_for_para_problem}
    Consider problem (\ref{eqn:para_min_problem}). Assume:
        (a) 
        for each $(y,t)\in Y\times T$, 
        $X(y,t)\subset \R^n$ is a compact and convex set; and there exist a positive constant $A$ and a compact set 
        $\hat{X}$ such that
        $\inmat{diam}(X(y,t))\leq A$ and $\cup_{y\in Y,t\in T}X(y,t)\subset \hat{X}$;
        (b) the Slater's condition (\ref{eqn:slater})  holds 
        uniformly for all $(y,t)\in Y\times T$ 
        with a positive constant $\alpha$; 
        (c) for each  fixed $(y,t)$, 
         $\psi_i(x,y,t)$ is convex in $x$
        for $i=1,\dots,m$;
        (d)  $\psi(x,y,t)$ is
        H{\"o}lder continuous in $t$ 
        uniformly for $(x,y)\in \hat{X}\times Y$ {\color{black}:}
        there exist positive constants $\kappa_\psi$ and $q_\psi\in (0,1]$ such that
        \bgeq
            \left\| \psi(x,y,t')-\psi(x,y,t'') \right\| \leq \kappa_\psi \| t' - t''\| ^{q_\psi}, \forall t', t''\in T 
        \edeq
       uniformly for all $(x,y)\in \hat{X}\times Y$;
        (e) 
        $\phi(x, y,t)$ 
        satisfies the following condition over $\hat{X} \times Y\times T$ {\color{black}:} there are positive constants $\kappa_\phi$ and $q_1,q_2\in (0,1]$ such that for all $y\in Y$,
        \bgeq
            |\phi(x', y, t') - \phi(x'', y, t'')| \leq \kappa_\phi(\|x' - x''\|^{q_1} +  \|t' - t''\|^{q_2}), \; \forall x', x'' \in \hat{X}, t', t'' \in T.
        \edeq
    Then the following assertions hold.
    \begin{itemize}
        \item[(i)] If (f) $\phi(x,y,t)$ is strictly quasi-convex in $x$, then the optimal solution set $X^*(y,t)$ is a singleton $\{x^*(y,t)\}$ for each $(y,t)$, and the optimal value function
        $\phi^*(y,t)$ and the optimal solution mapping $x^*(y,t)$ are continuous in $t\in T$.
    \item[(ii)] If 
    (f') $\phi(x,y,t)$ satisfies the growth condition in $x$ {\color{black}:} there exists positive constants $\nu,\beta$ such that for all $(y,t)$,
    \bgeq
        \phi(x,y,t) \geq \phi^* (y,t) + \beta \dd(x,X^*(y,t))^\nu, \forall x\in X(y,t),
    \edeq
    then, for all $y\in Y$,
        \bgeqn
    \label{eqn:continuity_of_optimal_solution}
        &&\mathbb{H}(X^*(y,t'),X^*(y,t'')) \\\nonumber
        &\leq& \frac{A\kappa_\psi}{\alpha}\|t''-t'\|^{q_\psi}+ 
\left(\frac{4\kappa_\phi}{\beta}\left( \left(\frac{A\kappa_\psi}{\alpha}\right)^{q_1}+1\right)\right)^{1/\nu}
\|t'-t''\|^{q/\nu}\\\nonumber
&\leq&\left(\frac{A\kappa_\psi}{\alpha}+ 
\left(\frac{4\kappa_\phi}{\beta}\left( \left(\frac{A\kappa_\psi}{\alpha}\right)^{q_1}+1\right)\right)^{1/\nu}\right)
\|t'-t''\|^{\min \{q_\psi,q/\nu\}}
,\\\nonumber
&&\qquad\qquad\qquad\qquad\qquad\qquad\qquad\qquad\forall t',t''\in T, \;\inmat{with} \; \|t'-t''\|\leq 1,
    \edeqn
    where $q := \min\{q_1q_\psi,q_2\}$.
    \end{itemize}
\end{theorem}

\noindent
\textbf{Proof.} 
{\color{black}
Part (i). 
Under conditions (a)-(c), 
we know from Lemma \ref{lemma:hoffman} (ii)
that the feasible set-valued mapping $X(y,t)$ is 
continuous in $t\in T$ 
uniformly for all $y\in Y$.
Condition (e) states that the objective function is continuous in $(x,t)$. 
By Berge's maximum theorem \cite{berge1963topological}, the optimal value 
function $\phi^*(y,t)$ is continuous in $t\in T$, 
and the set
of optimal solutions $X^*(y,t)$ is upper semi-continuous in $t\in T$. 
Moreover, under conditions (c) and 
(f), $X^*(y,t)=\{x^*(y,t)\}$ is a singleton for each fixed $(y,t)$, and thus optimal solution mapping $x^*(y,t)$ is continuous in $t$. 

Part (ii). For any $y\in Y$, let $t',t''\in T$.
Let $x^*(y,t')\in X^*(y,t')$, $x^*(y,t'') \in X^*(y,t'')$,
$x' \in \arg\min\limits_{x\in X(y,t'')} \|x-x^*(y,t')\|$.
Then
\bgeqn
\label{eqn:proof_part2_equicont}
    && \phi(x^*(y,t''),y,t'')- \phi(x^*(y,t'),y,t')\\
    & \leq & \phi(x',y,t'')- \phi(x^*(y,t'),y,t')\nonumber\\
    & \leq & \kappa_\phi(\|x' - x^*(y,t')\|^{q_1} +  \|t'' - t'\|^{q_2})\nonumber\\
    & \leq & \kappa_\phi(\mathbb{H}(X(y,t'),X(y,t'')) ^{q_1} +  \|t' - t''\|^{q_2})\nonumber\\
    & \leq & \kappa_\phi\left( \left(\frac{A\kappa_\psi}{\alpha}\right)^{q_1}\| t' - t''\|^{q_\psi q_1} +  \|t' - t''\|^{q_2}\right)\nonumber\\
     & \leq & \kappa_\phi\left( \left(\frac{A\kappa_\psi}{\alpha}\right)^{q_1}+1\right) \| t' - t''\|^q,\nonumber
\edeqn
where 
the second inequality 
 {\color{black}arises} from the condition (e), 
the third 
 {\color{black}arises}
from the definition of $x'$: $\|x' - x^*(y,t')\| = \dd(x^*(y,t'),X(y,t'')) \leq \mathbb{D}(X(y,t'),X(y,t''))$, 
the fourth 
 {\color{black}arises}
from (\ref{eqn:stability_of_feasible_set}) in Lemma \ref{lemma:hoffman} (ii), and the fifth is obtained by letting $q := \min\{q_1q_\psi,q_2\}\in (0,1]$,
and 
$\| t''-t' \|\leq \delta \leq 1$.

Let $x''\in X(y,t'')\subset X(y,t') + \frac{A\kappa_\psi}{\alpha}\|t'-t''\|^{q_\psi}\B$, 
where $\B$ denotes the unit ball in $\R^n$. 
By the definition of $q$, we have $q_{\psi} \geq \frac{q}{q_1}$ and thus $x'' \in X(y,t') + \frac{A\kappa_\psi}{\alpha}\|t'-t''\|^{q/q_1}\B$.
Let $\epsilon$ be a positive number to be specified later.
Suppose that $\dd(x'',X^*(y,t'))\geq \epsilon$.
We consider two cases. 

Case 1: $x''\in X(y,t')$. By condition (f'), we have
\bgeqn 
\phi(x'',y,t') -
\phi(x^*(y,t'),y,t') \geq \beta\epsilon^\nu.
\label{eqn:proof-case1}
\edeqn
Using condition (e), inequalities
(\ref{eqn:proof_part2_equicont}) and (\ref{eqn:proof-case1}), we obtain
\bgeqn
\label{eq:Lema-BNBE-proof-case1-a}
&&\phi(x'',y,t'') -
\phi(x^*(y,t''),y,t'') \\\nonumber
&=& 
\phi(x'',y,t') -
\phi(x^*(y,t'),y,t')
+(
\phi(x'',y,t'') -
\phi(x'',y,t')
)\\
&&- (\phi(x^*(y,t''),y,t'')-\phi(x^*(y,t'),y,t'))
\nonumber\\
&\geq& \beta\epsilon^\nu -
\kappa_\phi \|t'-t''\|^{q_2} -
\kappa_\phi\left( \left(\frac{A\kappa_\psi}{\alpha}\right)^{q_1}+1\right) \| t' - t''\|^q\nonumber\\
&\geq& \beta\epsilon^\nu - \kappa_\phi\left( \left(\frac{A\kappa_\psi}{\alpha}\right)^{q_1}+2\right)\|t'-t''\|^q. \nonumber
\edeqn 

Case 2: $x''\not\in X(y,t')$. Let
$\hat{x}'' \in \arg\min\limits_{x\in X(y,t')} \|x-x''\|$,
which means $\|\hat{x}''-x''\| = \dd(x'',X(y,t'))\leq \frac{A\kappa_\psi}{\alpha}\|t'-t''\|^{q_\psi} \leq  \frac{A\kappa_\psi}{\alpha}\|t'-t''\|^{q/q_1}$,
where the inequality 
 {\color{black}arises}
from the definition of $x''$.
Thus
\bgeqn
\label{eq:Lema-BNBE-proof-case2-a}
&&\phi(x'',y,t'') -
\phi(x^*(y,t''),y,t'') \\
&=& 
\phi(\hat{x}'',y,t') -
\phi(x^*(y,t'),y,t') +
\phi(x'',y,t') -\phi(\hat{x}'',y,t')
\nonumber\\
&& 
+(
\phi(x'',y,t'') -
\phi(x'',y,t')
)
- (\phi(x^*(y,t''),y,t'')-\phi(x^*(y,t'),y,t'))\nonumber\\
&=& R + 
\phi(x'',y,t') -\phi(\hat{x}'',y,t'), \nonumber
\edeqn
where 
\bgeq 
R &:=& \phi(\hat{x}'',y,t') -
\phi(x^*(y,t'),y,t') \\
&&
+(
\phi(x'',y,t'') -
\phi(x'',y,t')
)
- (\phi(x^*(y,t''),y,t'')-\phi(x^*(y,t'),y,t')). 
\edeq
As in Case 1, we have by virtue of (\ref{eq:Lema-BNBE-proof-case1-a}), 
\bgeqn 
R &\geq& \beta(\dd(\hat{x}'',X^*(y,t'))^\nu - \kappa_\phi\left( \left(\frac{A\kappa_\psi}{\alpha}\right)^{q_1}+2\right)\|t'-t''\|^q\nonumber\\
&\geq& 
\beta(\dd(x'',X^*(y,t'))-\|\hat{x}''-x''\|)^\nu - \kappa_\phi\left( \left(\frac{A\kappa_\psi}{\alpha}\right)^{q_1}+2\right)\|t'-t''\|^q\nonumber\\
&\geq& \beta\left(\epsilon-
\frac{A\kappa_\psi}{\alpha}\|t''-t'\|^{q_\psi}\right)^\nu - \kappa_\phi\left( \left(\frac{A\kappa_\psi}{\alpha}\right)^{q_1}+2\right)\|t'-t''\|^q.
\label{eq:Lema-BNBE-proof-case2-b}
\edeqn 
On the other hand, by condition (e),
\bgeqn 
|\phi(x'',y,t') -\phi(\hat{x}'',y,t')|
&\leq& 
\kappa_\phi\|x''-\hat{x}''\|^{q_1}
\leq \kappa_\phi\left(\frac{A\kappa_\psi}{\alpha}\right)^{q_1}\|t''-t'\|^{q}.
\label{eq:Lema-BNBE-proof-case2-c}
\edeqn 
Combining (\ref{eq:Lema-BNBE-proof-case2-a})
-(\ref{eq:Lema-BNBE-proof-case2-c}), we obtain
\bgeq
&&\phi(x'',y,t'') -
\phi(x^*(y,t''),y,t'') \\
&\geq& 
\beta\left(\epsilon-
\frac{A\kappa_\psi}{\alpha}\|t''-t'\|^{q_\psi}\right)^\nu - 2\kappa_\phi\left( \left(\frac{A\kappa_\psi}{\alpha}\right)^{q_1}+1\right)\|t'-t''\|^q.
\edeq

\begin{figure}[t]
    \centering
    \includegraphics[width=0.4\textwidth]{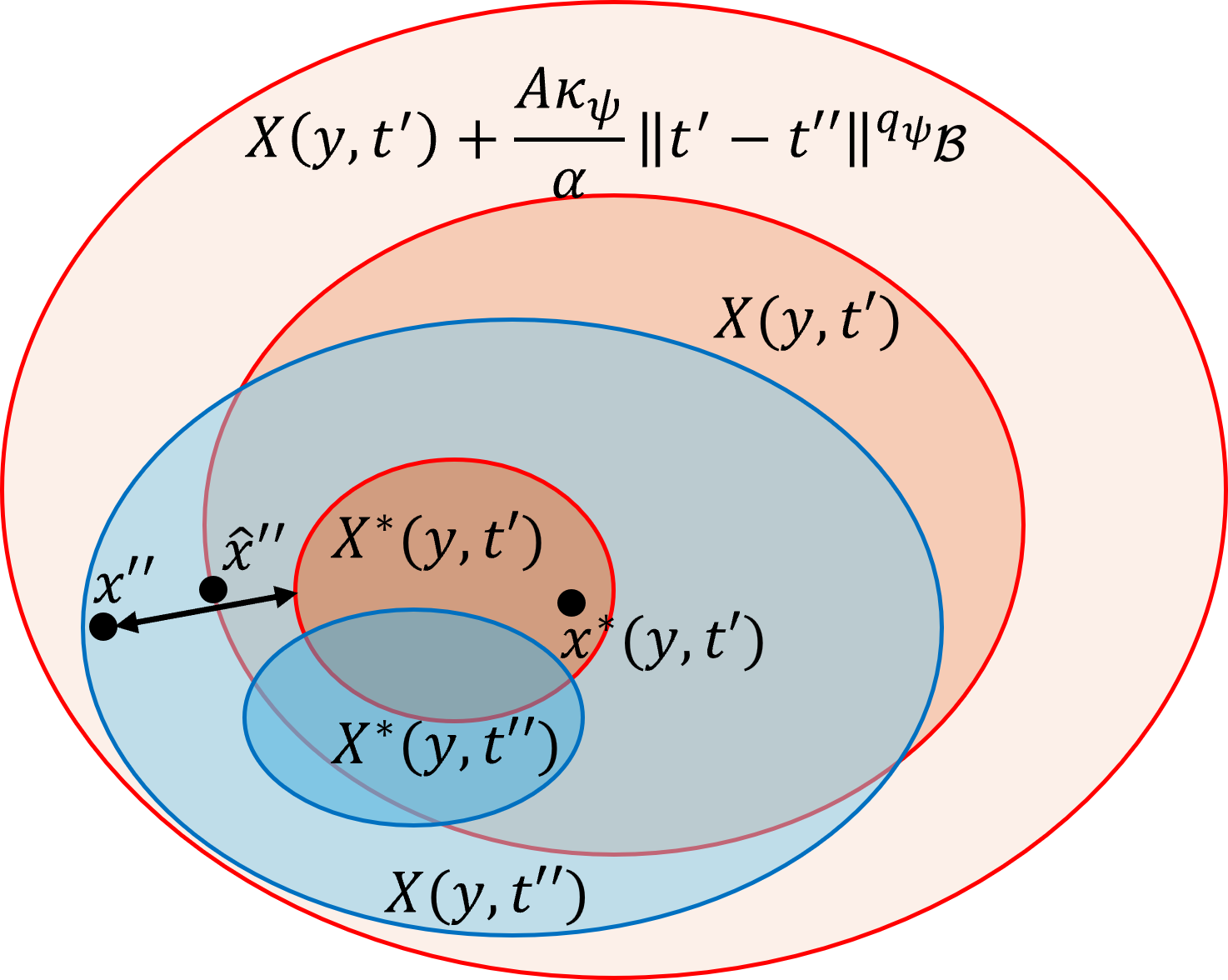}
    \caption{
    \small Illustration of the proof of Case 2 in Part (ii), where notation `$\longleftrightarrow$' denotes $\dd(x'',X^*(y,t'))$ which is larger than $\epsilon$.}
    \label{fig:continuity}
    \vspace{-15pt}
\end{figure}

Summarizing the discussions above, we 
may set
\bgeqn 
\epsilon = \frac{A\kappa_\psi}{\alpha}\|t''-t'\|^{q_\psi}+ 
\left(\frac{4\kappa_\phi\left( \left(\frac{A\kappa_\psi}{\alpha}\right)^{q_1}+1\right)\|t'-t''\|^q}{\beta}\right)^{1/\nu}
\edeqn 
and 
then obtain
$
\phi(x'',y,t'') -
\phi(x^*(y,t''),y,t'') >0,
$
which means $x''\not\in X^*(y,t'')$. Since
$X^*(y,t'') \subset X(y,t'')\subset X(y,t') + \frac{A\kappa_\psi}{\alpha}\|t'-t''\|^{q_\psi}\B$, the inequalities above imply that 
for any $x^*(y,t'')\in X^*(y,t'')$,
we must have $\dd(x^*(y,t''),X^*(y,t')) < \epsilon$, and thus
$\mathbb{D}(X^*(y,t''),X^*(y,t')) < \epsilon$.
By swapping the position between $X^*(y,t')$ and $X^*(y,t'')$, we immediately arrive at (\ref{eqn:continuity_of_optimal_solution})
\hfill $\Box$
}

In the case 
that $\psi(x,y,t)$ and $\phi(x,y,t)$ are (uniformly) Lipschitz continuous reduced from H{\"o}lder continuous {\color{black}:} $q_\psi,q_1,q_2=1$, 
the parameters in Theorem~\ref{thm:stability_for_para_problem} 
can be simplified with
$q_\psi,q_1,q=1$.
By Theorem \ref{thm:stability_for_para_problem} (ii), we can 
derive the equicontinuity results (\ref{eqn:equicontinuity}). 
We are now ready to establish the equicontinuity of the optimal response functions of players in  {\color{black}the} GBNE model and 
 {\color{black}also}
the existence of the GBNE.



\subsection{Existence of a continuous generalized Bayesian Nash equilibrium}
With the technical preparations in  section \ref{subsec:continuity_of_para_program}, we turn to discuss the existence of a continuous GBNE of model (\ref{eqn:GNE}).
We do so by virtue of Schauder's fixed point theorem 
as in \cite{guo2021existence} and \cite{meirowitz2003existence}. 
To this end, we recall
the basics related to the classic result.

Let $Z$ be a Banach space and $T:Z\to Z$ be an operator.
A set in $Z$ is called {\em relatively compact} if its closure is compact. 
The operator $T$ is said to be {\em compact}
if it is continuous and maps the bounded sets into the relatively compact sets.

\begin{theorem}[Schauder's fixed point theorem \cite{zeilder1985nonlinear}]
If $\mathcal{Z}$ is a non-empty, closed, bounded, convex subset of a Banach space $Z$
and $T: \mathcal{Z}\to \mathcal{Z}$ is a compact operator,
then $T$ has a fixed point.
\label{t-Shauder}
\end{theorem}

It is well-known that the relative compactness 
may be characterized 
by the uniform boundedness and equicontinuity. The next proposition states this.

\begin{proposition} [Arzela-Ascoli Theorem 
 {\color{black}Theorem A5 in \cite{rudin1973functional}}]
\label{lemma:arzela-ascoli}
    A set $\mathcal{Z}\subset Z$ is relatively compact if and only if
$\mathcal{Z}$ satisfies  the following two conditions:
\begin{itemize}
    \item[(a)] uniform boundedness {\color{black}:}
$
\sup_{f\in \mathcal{Z}}
 \|f\|_\infty
 <\infty,
$
and
\item[(b)] equicontinuity {\color{black}:}
for every $\epsilon>0$, there exists a constant $\delta>0$ such that
$
\sup_{f\in \mathcal{Z}} \|f(\theta)-f(\theta')\| < \epsilon, \forall \theta, \theta' \; \inmat{with}\; \|\theta-\theta'\|<\delta.
$
\end{itemize}

\end{proposition}

 {\color{black}Here, we discuss the compactness/relative compactness of a set in the context of the norm-induced topology.
In general, the strong compactness 
in Banach spaces is a stringent requirement.
In this context, we avoid directly assuming the compactness or relative compactness in Banach spaces. 
Instead, we derive the relative compactness of the set of the player's response strategies 
by virtue of Arzela-Ascoli Theorem.}

 {\color{black}
Before proceeding to the discussion on existence of GBNE, 
we recall the 
model and 
associated 
notation
to facilitate reading.  
By letting
\bgeq
    \rho_i(a_i,f_{-i},\theta_i)& :=& \int_{\Theta_{-i}}u_i(a_i,f_{-i}(\theta_{-i}),\theta )d\eta_i(\theta_{-i}|\theta_i),\\
    h_i(a_i,f_{-i},\theta_i)& :=& \int_{\Theta_{-i}}g_i(a_i,f_{-i}(\theta_{-i}),\theta) d\eta_i(\theta_{-i}|\theta_i),
\edeq
we can write the 
GBNE model as
\bgeq
        f_i^*(\theta_i)\in\arg \max\limits_{a_i\in \mathbb{A}_i} \; \rho_i(a_i,f_{-i}^*,\theta_i), 
        \inmat{s.t.} \; h_i(a_i,f_{-i}^*,\theta_i)\leq 0, 
\edeq
for each $i\in N$, $\theta_i\in \Theta_i$, see
(\ref{eqn:GNE}).
Moreover, 
we define the mapping:
\bgeq
 \Psi (f)&:=& \bigg\{ (\tilde{f}_1,\dots,\tilde{f}_n)\in \mathcal{F}:
    \inmat{for} \; i\in N, 
    \\
    && \tilde{f}_i(\theta_i)\in \arg \max_{a_i\in \mathcal{A}_i(f_{-i},\theta_i)} \rho_i(a_i,f_{-i},\theta_i), \quad \forall \theta_i\in \Theta_i\bigg\},
\edeq
where $\mathcal{A}_i(f_{-i},\theta_i) := \left\{ a_i\in \mathbb{A}_i: \int_{\Theta_{-i}}g_i(a_i,f_{-i}(\theta_{-i}),\theta)d\eta_i(\theta_{-i}|\theta_i)\leq 0 \right\}$.
Then $f^*$ is a GBNE iff 
$f^* \in \Psi (f^*)$.
}

In the forthcoming discussions, we show that $\Psi$ defined in (\ref{eqn:mapping})
has a fixed point when it is applied to a subset of 
${\cal F}$,
 {\color{black} and thus the corresponding generalized Bayesian game has an GBNE.}
For this purpose, we derive sufficient conditions
for the continuity of $\Psi$ and relative compactness of $\Psi(\cal F)$. 
By Proposition \ref{lemma:arzela-ascoli}, we provide conditions which guarantee
the functions in set $\Psi({\cal F})$ are uniformly bounded and equicontinuous.
The next assumption prepares for this.


\begin{assumption}\label{assump:GNE}
    Consider the model (\ref{eqn:GNE}). For $i\in N$
    \begin{itemize}
        \item[(a)] $\Theta_i$ and $\mathbb{A}_i$ are convex and compact;
        \item[(b)] $\rho_i(a_i,f_{-i},\theta_i)$ is continuous in $(a_i,f_{-i},\theta_i)$, and $h_i(a_i,f_{-i},\theta_i)$ is continuous in $(a_i,f_{-i},\theta_i)$;
        \item[(c)] for each $f_{-i}\in \mathcal{F}_{-i}$ and $\theta_i\in \Theta_i$, $\rho_i(a_i,f_{-i},\theta_i)$ is strictly quasi-concave on $\mathcal{A}_{i}(f_{-i},\theta_i)$, and $h_i(a_i,f_{-i},\theta_i)$ is convex in $a_i$; 
       \item[(d)] there exist positive constants $\kappa_\rho,\kappa_h,q_1, q_2,q_h$ such that for all $f_{-i}:\Theta_{-i}\rightarrow \mathbb{A}_{-i}$,
        \bgeq
            | \rho_i(a_i,f_{-i},\theta_i) - \rho_i(a_i',f_{-i},\theta_i^\prime)| &\leq&  \kappa_\rho \left(\|a_i - a_i'\|^{q_1} + \|\theta_i - \theta_i^\prime\|^{q_2}\right),\\
            &&\qquad\qquad\qquad\;\forall a_i,a_i'\in \mathbb{A}_i, \theta_i,\theta_i'\in \Theta_i,\\
            \| h_i(a_i,f_{-i},\theta_i)-h_i(a_i,f_{-i},\theta_i') \| &\leq& \kappa_h \| \theta_i - \theta_i^\prime\|^{q_h},
            \forall a_i,a_i'\in \mathbb{A}_i, \theta_i,\theta_i'\in \Theta_i;
        \edeq
        \item[(e)] 
        for any $f_{-i} \in \mathcal{F}_{-i}, \theta_i\in \Theta_i$, the feasible set $\mathcal{A}_i(f_{-i},\theta_i)$ is a non-empty and compact set with 
        $\inmat{diam}(\mathcal{A}_i(f_{-i},\theta_i))< \inmat{diam}(\mathbb{A}_i)\leq A,$  {\color{black} where $A:=\max\limits_{i\in N} \inmat{diam}(\mathbb{A}_i)$}.
        Moreover, 
        Slater's condition holds for all parametric inequality systems in model (\ref{eqn:GNE}) {\color{black}: there exists a positive constant $\alpha>0$ such that for any $(f_{-i},\theta_i)$, there exists a feasible solution, denoted by $\bar{a}_i(f_{-i},\theta_i)\in \mathbb{A}_i$, satisfying\linebreak
        $h_i(\bar{a}_i(f_{-i},\theta_i),f_{-i},\theta_i)\leq -\alpha$};
        
        \item[(f)] $\rho_i(a_i,f_{-i},\theta_i)$ satisfies the  {\color{black} following} growth condition in $a_i$:
         {\color{black} there exist positive parameters $\beta, \nu$ such that }
        \bgeq
            -\rho_i(a_i^\prime,f_{-i},\theta_i) \geq - \rho_i^*(f_{-i},\theta_i)+\beta d(a_i^\prime, \mathcal{A}_i^*(f_{-i},\theta_i))^\nu, \forall a_i^\prime \in \mathbb{A}_i,
        \edeq
        uniformly for all $f_{-i}\in \mathcal{F}_{-i}$ and $\theta_i\in \Theta_i$,  where $\rho_i^*(f_{-i},\theta_i)$ and $\mathcal{A}_i^*(f_{-i},\theta_i)$ denote the optimal value and optimal solution set respectively:
        \bgeq
            \rho^*_i(f_{-i},\theta_i) &:= & \max_{a_i\in \mathcal{A}_i(f_{-i},\theta_i)} \rho_i(a_i,f_{-i},\theta_i),\\
            \mathcal{A}_i^*(f_{-i},\theta_i) &:=& \arg\max_{a_i\in \mathcal{A}_i(f_{-i},\theta_i)} \rho_i(a_i,f_{-i},\theta_i).
        \edeq
    \end{itemize}
\end{assumption}


\begin{remark}
Conditions (a) and (f) are used in \cite{guo2021existence},
conditions (c)-(e) are newly introduced to address the dependence of 
feasible sets $\mathcal{A}_i(f_{-i},\theta_i)$ on $(f_{-i},\theta_i)$ for $i\in N$. We make 
 {\color{black}the following remarks about the latter.}
\begin{enumerate}
        \item[(i)] Condition (c) ensures that 
        the feasible set $\mathcal{A}_i(f_{-i},\theta_i)$ is a convex set and
        the set of optimal solutions 
        $\mathcal{A}_i^*(f_{-i},\theta_i)$ is a singleton. 
        The quasi-concavity condition of $\rho_i(\cdot,f_{-i},\theta_i)$ over $\mathcal{A}_i(f_{-i},\theta_i)$
        may be verified by the quasi-concavity of the function 
        over 
        set $\mathbb{A}_i$ which contains $\mathcal{A}_i(f_{-i},\theta_i)$.
         
        \item[(ii)] Condition (d) 
        requires
        uniform H{\"o}lder continuity of the
        objective 
        function 
        and
        constraint function. 
        The former extends the Lipschtz continuity in \cite{guo2021existence}, the latter is 
        new.
This condition 
enables us to use Theorem \ref{thm:stability_for_para_problem} to derive the uniform continuity of $\mathcal{A}_i^*(f_{-i},\theta_i)$.

        \item[(iii)] Condition (e) ensures the non-emptiness
        of $\mathcal{A}_i(f_{-i},\theta_i)$,
       see the similar 
       condition 
       in \cite{ichiishi1983game} for generalized Nash games. The assumption  {\color{black}about} the compactness of each player's feasible set can be weakened 
       to the super compactness condition as
       in \cite{guo2021existence}, 
       here we 
       make it directly for 
       simplicity of presentation. 
        Slater's condition is widely used in the literature of optimization. 
    \end{enumerate}
Note that it is possible to present 
stronger but more verifiable 
sufficient conditions for (b), (c), (d), (f), e.g.,
by considering that $u_i(a_i,a_{-i},\theta)$ is strongly concave in $a_i$ for each $\theta\in\Theta$, see 
\cite{guo2021existence,ravat2011characterization,ravat2017existence}.    
\end{remark}
   
We are now ready to state the main result of this section.

\begin{theorem}[Existence of continuous GBNE] \label{thm:existence} 
Consider the GBNE model (\ref{eqn:GNE}). Under Assumption \ref{assump:GNE}, the 
problem has a generalized Bayesian Nash equilibrium with the best strategy functions being uniformly H\"older continuous.

\end{theorem}

\noindent
\textbf{Proof.} 
As explained  {\color{black}above,}
it suffices to show
that the mapping $\Psi$ defined
 in (\ref{eqn:mapping}) has 
 a fixed point which is a GBNE. For $i\in N$, let 
 $\tilde{\mathcal{C}}_i$ be a set of continuous functions 
 $\tilde{f}_i:\Theta_i\to \mathcal{A}_i(f_{-i},\theta_i)$, 
and
 $\tilde{\mathcal{C}} :=\prod_{i\in N}\tilde{\mathcal{C}}_i$. 
 Under Assumption \ref{assump:GNE} (a) and (e), $\tilde{\mathcal{C}}$ is a non-empty, closed, bounded and convex subset of a Banach space
equipped with the infinity norm. We proceed  {\color{black}with} the rest of the proof in three steps. 

\underline{Step 1.}  We prove that $\mathcal{A}_i^*(f_{-i},\theta_i)$ is continuous in $(f_{-i},\theta_i)$. 
Observe first that 
condition (e) ensures the non-emptiness and compactness of 
the feasible set $\mathcal{A}_i(f_{-i},\theta_i)$ 
whereas the convexity of $h_i(a_i,f_{-i},\theta_i)$ with respect to  $a_i$  {\color{black} which follows from} condition (c) 
ensures the convexity of the feasible set.
Moreover, the strict quasi-concavity of $\rho_i(a_i,f_{-i},\theta_i)$ in $a_i$ 
 {\color{black} according to} condition (c) 
guarantees that  $\mathcal{A}_i^*(f_{-i},\theta_i)$ is a singleton.
On the other hand,
under the Slater's condition and the uniform H\"older continuity of $h_i(a_i,f_{-i},\theta_i)$ in $\theta_i$,
we know by virtue of Lemma \ref{eqn:para_min_problem} (ii)
that the feasible set $\mathcal{A}_i(f_{-i},\theta_i)$ is continuous in $(f_{-i},\theta_i)$.
By the classic stability results (e.g. Theorem in 4.2.1 \cite{bank1984Nonlinear} or Theorem 2.3.1 in \cite{ichiishi1983game}), we conclude that  $\mathcal{A}_i^*(f_{-i},\theta_i)$
is continuous in $(f_{-i},\theta_i)$.

\underline{Step 2.} We show that $\mathcal{A}_i^*(f_{-i},\theta_i)$ is equicontinuous over $\Theta_i$ for all $f_{-i}\in \tilde{\mathcal{C}}_i$, and the set $\Psi(\tilde{C})$ is relatively compact.
Under 
 {\color{black}Assumption \ref{assump:GNE}},
we can use
Theorem \ref{thm:stability_for_para_problem} to assert
 that the optimal solution mapping
 $\mathcal{A}_i^*(f_{-i},\theta_i)$ is equicontinuous over $\Theta_i$ for all $f_{-i}\in \tilde{\mathcal{C}}_i$ {\color{black}:}
\bgeqn
\label{eqn:equicontinuity_of_optimal_solution_set}
     \|\mathcal{A}_i^*(f_{-i},\theta_i) - \mathcal{A}_i^*(f_{-i},\theta_i^\prime)\| 
    \leq \kappa \|\theta_i-\theta_i'\|^q
    , \forall \theta_i,\theta_i^\prime \in \Theta_i,  {\color{black}\; \inmat{with} \; \|\theta_i-\theta_i'\|\leq 1, }
\edeqn
where $q = \min\{q_h,\frac{q_1q_h}{\nu},\frac{q_2}{\nu}\}$, 
 {\color{black}
$
\kappa:=
\left(\frac{A\kappa_h}{\alpha}+ 
\left(\frac{4\kappa_\rho}{\beta}\left( \left(\frac{A\kappa_h}{\alpha}\right)^{q_1}+1\right)\right)^{1/\nu}\right).
$}
Moreover, 
since 
condition (a), (e) guarantees that $\cup_{\theta_i\in\Theta_i}\mathcal{A}_i(f_{-i},\theta_i)\subset \mathbb{A}_i$ is uniformly bounded,
then
$\cup_{\theta_i\in\Theta_i}\mathcal{A}_i^*(f_{-i},\theta_i) \subset \cup_{\theta_i\in\Theta_i}\mathcal{A}_i(f_{-i},\theta_i)$ is also uniformly bounded.
 By Proposition \ref{lemma:arzela-ascoli}, 
 $\Psi(\tilde{C})$ is relatively compact.
\underline{Step 3.}
Based on the continuity of 
$\mathcal{A}_i^*(f_{-i},\theta_i)$ in $(f_{-i},\theta_i)$
and  {\color{black}the} equicontinuity of
$\mathcal{A}_i^*(f_{-i},\theta_i)$ 
over $\Theta_i$, 
we show that $\Psi:\tilde{C}\rightarrow \tilde{C}$ is a continuous operator. 
Since $\Theta_i$ is compact under condition (a), there is a $\delta$-net covering $\Theta_i$: for any small positive number  {\color{black}$\delta\in (0,1)$}, there exists a finite number of points $\theta_i^1,\dots,\theta_i^M \in \Theta_i$ such that for any $\theta_i\in\Theta_i$, there is $k\in\{1,\dots,M\}$, such that $\|\theta_i-\theta_i^k\|\leq \delta$.
Moreover, by the continuity of $\mathcal{A}_i^*(f_{-i},\theta_i)$ in $(f_{-i},\theta_i)$ proved in Step 1, 
for any $\epsilon>0$ and fixed $f_{-i}\in \tilde{\mathcal{C}}_i$,  
we have $\|\mathcal{A}_i^*(f_{-i},\theta_i^k)-\mathcal{A}_i^*(\tilde{f}_{-i},\theta_i^k)\| \leq \epsilon$, for $k= 1,\dots,M$ when
 $\tilde{f}_{-i}$ is sufficiently close to $f_{-i}$.
Together with 
(\ref{eqn:equicontinuity_of_optimal_solution_set}),
we have 
\bgeq
    &&\|\mathcal{A}_i^*(\tilde{f}_{-i},\theta_i)-\mathcal{A}_i^*(f_{-i},\theta_i)\|\\
    & \leq & 
    \|\mathcal{A}_i^*(\tilde{f}_{-i},\theta_i)-\mathcal{A}_i^*(\tilde{f}_{-i},\theta_i^k)\|
    + \|\mathcal{A}_i^*(\tilde{f}_{-i},\theta_i^k)-\mathcal{A}_i^*(f_{-i},\theta_i^k)\|\\
    && +
    \|\mathcal{A}_i^*(f_{-i},\theta_i^k)-\mathcal{A}_i^*(f_{-i},\theta_i)\|  
 \leq  2\kappa \delta^q + \epsilon, 
\edeq
and hence 
$\sup_{\theta_i\in\Theta_i}\|\mathcal{A}_i^*(\tilde{f}_{-i},\theta_i)- \mathcal{A}_i^*(f_{-i},\theta_i)\|\leq 2\epsilon$ for $\delta \leq \left( \frac{\epsilon}{2\kappa} \right)^{\frac{1}{q}}$. 
Thus, 
$\mathcal{A}_i^*(\cdot,\theta_i):
\tilde{C}_{-i}\rightarrow \tilde{C}_i$, $i\in N$, is continuous and $\Psi:\tilde{C}\rightarrow \tilde{C}$ is a continuous operator. 
The relative compactness of $\Psi(\tilde{C})$ 
and continuity of $\Psi:\tilde{C}\rightarrow \tilde{C}$ 
allow us to claim by Theorem \ref{t-Shauder} that 
$\Psi:\tilde{C}\rightarrow \tilde{C}$ has a fixed point which is a GBNE, and the uniform H\"older continuity follows from (\ref{eqn:equicontinuity_of_optimal_solution_set}).
\hfill $\Box$

\subsection{Uniqueness of 
continuous GBNE} 
It is challenging to 
the
derive  {\color{black}the} uniqueness of a continuous GBNE 
 {\color{black}because} even in generalized Nash games,  {\color{black}the} uniqueness of (global) GNE can only be 
established under 
specific
circumstances, 
see e.g.~\cite{facchinei201012nash}. 
Here, we 
first 
consider  
the GBNE model where 
each player's feasible set of actions is only 
 type-dependent:
\begin{subequations}
\label{eqn:GBNE_type_dependent}
    \bgeqn
    f^*\in \arg & \max_{f\in \mathcal{C}} &\sum_{i=1}^n \mathbb{E}_{\eta} [u_i(f_i(\theta_i),f_{-i}^*(\theta_{-i}),
     {\color{black}\theta}
    )]\\
    & \inmat{s.t.} &  {\color{black}h}_i(f_i(\theta_i),\theta_i)\leq 0,\quad \forall \theta_i\in\Theta_i, i\in N.\label{eqn:GBNE_type_dependent_b}
\edeqn
\end{subequations}
By slightly abusing  {\color{black}the} notation, let  
$\mathcal{A}_i(\theta_i) := \{a_i\in \mathbb{A}_i: g_i(a_i,\theta_i)\leq 0\}$, $\tilde{\mathcal{C}}_i'$ be a set of continuous functions $f_i:\Theta_i\rightarrow \mathcal{A}_i(\theta_i)$, and $\tilde{\mathcal{C}}':= \prod_{i\in N}\tilde{\mathcal{C}}_i'$. 
 {\color{black}
There are two important cases for the constraints (\ref{eqn:GBNE_type_dependent_b}). 

Case 1. 
$g_i(a_i,\theta_i)$ only depends on player $i$'s own action and type.
The constraint reduces to 
    \begin{eqnarray}
        h_i(a_i,\theta_i) = g_i(a_i,\theta_i)\leq 0.
    \end{eqnarray}
The resulting model subsumes the BNBE model in Liu et al.~\cite{liu2024Bayesian}. 

Case 2. 
$g_i$ 
depends on $\theta_{-i}$, that is, the rivals' types.
 In that case, (\ref{eqn:GBNE_type_dependent_b}) can be formulated as 
 \begin{eqnarray} \label{def:constraint_explanation2}
h_i(f_i(\theta_i),\theta_i) = \int_{\Theta_{-i}}g_i(f_i(\theta_i),\theta)d\eta_i(\theta_{-i}|\theta_i) \leq 0,
\end{eqnarray}
where the expectation with respect to
 conditional distribution $\eta_{i}(\theta_{-i}|\theta_i)$ is needed because player $i$ 
 does not have complete information 
 for $\theta_{-i}$.
}



For the uniqueness of GBNE, We need the following condition. 


\begin{assumption} [Strict monotonicity] \label{assume:strictly_monotone}
        For each fixed $(a_{-i},\theta)$, 
        $u_i(a_i,a_{-i},\theta)$ is 
        directionally differentiable in $a_i$, 
        and
        \bgeqn
            \int_\Theta \sum_{i=1}^n &&\left[(u_i)'_{a_i}(f_i(\theta_i),f_{-i}(\theta_{-i}), \theta; \tilde{f}_i(\theta_i)-f_i(\theta_i))\right.\nonumber\\
            && +\left. (u_i)'_{a_i}(\tilde{f}_i(\theta_i),\tilde{f}_{-i}(\theta_{-i}), \theta; f_i(\theta_i)-\tilde{f}_i(\theta_i)) \right] \eta(d \theta)>0, \forall f,\tilde{f}\in \tilde{\mathcal{C}}', f \neq \tilde{f},
            \label{assumption:monotone_nondiff}
        \edeqn 
        where $(u_i)'_{a_i}(t,a_{-i},\theta;h)$ denotes the directional derivative with respect to $a_i$ along direction $h$ at a point $a_i=t$. 
\end{assumption}

In the case that
$u_i(a_i,f_{-i},\theta_i)$ is continuously differentiable with respect to $a_i$, 
 {\color{black} let $\nabla u (a,a,\theta):=(\nabla_{a_1} u_1(a_1,a_{-1},\theta), \dots, \nabla_{a_n} u_n(a_n,a_{-n},\theta))^T$. Then}
(\ref{assumption:monotone_nondiff}) reduces to
\bgeqn
    \;\qquad\mathbb{E}_\eta[ (\nabla u(f(\theta),f(\theta),\theta)-\nabla u(\tilde{f}(\theta),\tilde{f}(\theta),\theta)) ^T (f(\theta) - \tilde{f}(\theta))] < 0,
    \forall 
    f,\tilde{f}\in \tilde{\mathcal{C}}, f \neq \tilde{f}.\label{assumption:monotone_diff}
\edeqn
Moreover, if in addition $\nabla_{a_i} u_i (a_i,a_{-i},\theta)$ is continuously differentiable in $(a_i,a_{-i})$  a.e.~$\theta\in \Theta$ for all $i\in N$, then condition (\ref{assumption:monotone_diff}) is guaranteed by the negative 
definiteness of the Jacobian matrix of 
$\nabla u (a_i,a_{-i},\theta)$, see Lemma 4 in \cite{ui2016bayesian}.

 {\color{black}Given} this assumption, we are ready 
to derive the uniqueness of a continuous GBNE for problem (\ref{eqn:GBNE_type_dependent}).
\begin{theorem}[Uniqueness of the continuous GBNE in the type-dependent feasible set] \label{thm:uniqueness_type_dependent} 
Consider problem (\ref{eqn:GBNE_type_dependent}).
Under Assumption \ref{assume:strictly_monotone}, problem (\ref{eqn:GBNE_type_dependent}) has at most one continuous GBNE.
\end{theorem}

\noindent
\textbf{Proof.} 
Suppose for the sake of contradiction that there are two
distinct GBNEs, denoted by
$f,\tilde f\in \tilde{\mathcal{C}}$. 
Then
\bgeqn
\label{eqn:first-order}
    &&\int_\Theta  \sum_{i=1}^n (u_i)'_{a_i}(f_i(\theta_i),f_{-i}(\theta_{-i}), \theta;  \tilde{f}_i(\theta_i)-f_i(\theta_i)) \eta(d \theta)\\\nonumber
&=&  \sum_{i=1}^n \int_{\Theta_i} \left[\int_{\Theta_{-i}} (u_i)'_{a_i}(f_i(\theta_i),f_{-i}(\theta_{-i}), \theta;  \tilde{f}_i(\theta_i)-f_i(\theta_i))
d\eta_i(\theta_{-i}|\theta_i) \right] d\eta_i(\theta_i)\\\nonumber
&=&\sum_{i=1}^n \int_{\Theta_i} (\rho_i)'_{a_i}(f_i(\theta_i),f_{-i},\theta_i;  \tilde{f}_i(\theta_i)-f_i(\theta_i))d\eta_i(\theta_i)
\leq 0,
\edeqn
where the interchange of integration and directional differentiation  in the second 
equality follows from Page 79(4) in \cite{clarke1990optimization}, and the last inequality 
follows from
the first order optimality condition of $\rho_i$ at a GBNE $f_i$.
By swapping positions between $\tilde{f}$ and $f$, we have
\bgeqn 
   \int_\Theta  \sum_{i=1}^n (u_i)'_{a_i}(\tilde{f}_i(\theta_i),\tilde{f}_{-i}(\theta_{-i}), \theta; f_i(\theta_i)- \tilde{f}_i(\theta_i)) \eta(d \theta) \leq 0. 
   \label{eqn:first-order2}
\edeqn
Combining (\ref{eqn:first-order}) and (\ref{eqn:first-order2}), we immediately obtain a contradiction to (\ref{assumption:monotone_nondiff}).
\hfill $\Box$

Next, we discuss the case 
where each player's feasible set depends
not only on its own type
but also on its rivals' actions.
Under such circumstance, 
Assumption \ref{assume:strictly_monotone} is inadequate  
to guarantee
the uniqueness of GBNE. The next example explains this.

\begin{example} [Non-uniqueness of GBNE with the action-dependent feasible set]

    Consider a two-player generalized Bayesian game. Assuming  that type variables $\theta_1$ and $\theta_2$ are independent and uniformly distributed over $\Theta_1 =[1,2]$ and $\Theta_2 = [1,2]$, respectively,  {\color{black} we define a GBNE model by}
\bgeq
\begin{aligned}
    f_1^*(\theta_1) = \arg \max_{a_1\geq 0} &\; \int_0^1 \left[ 100 - a_1 - a_2(\theta_2) - \theta_1 \right] a_1 d\theta_2\\
    \inmat{s.t.} &\; \int_0^1 [ \theta_1 a_1 + \theta_2 a_2(\theta_2)  - 50 ] d\theta_2 \leq 0,\\
    &\; \theta_1 a_1 \leq 30,
\end{aligned}
\edeq
for $\theta_1\in\Theta_1= [1,2]$
and 
\bgeq
\begin{aligned}
    f_2^*(\theta_2) = \arg \max_{a_2\geq 0} &\; \int_0^1 \left[ 100 - a_1(\theta_1) - a_2 - \theta_2 \right] a_2 d\theta_1\\
    \inmat{s.t.} &\; \int_0^1 [ \theta_1 a_1(\theta_1) + \theta_2 a_2 - 50 ] d\theta_1 \leq 0,\\
    &\; \theta_2 a_2 \leq 30,
    \end{aligned}
\edeq
for $\theta_2\in\Theta_2= [1,2]$.
The Jacobian matrix of  $\nabla u (a_i,a_{-i},\theta)$ is 
$
\left[
\begin{array}{cc}
   -2  & -1 \\
   -1  & -2
\end{array}
\right],
$
which is negative definite and thus Assumption \ref{assume:strictly_monotone} is satisfied. 
However, we can verify that $f_1^*(\theta_1)= \frac{R}{\theta_1}, \theta_1\in [1,2]$, $f_2^*(\theta_2)= \frac{50-R}{\theta_2}, \theta_2\in [1,2]$ is a GBNE, for all $R \in [20,30]$, which means that 
the problem has infinitely many GBNEs. 
\end{example}

\section{Polynomial 
GBNE and SGNE}
\label{sec:polynomial}
Since 
the GBNE
is an infinite-dimensional 
equilibrium problem, 
it is usually difficult to compute it.
In this section, we propose a computational scheme which allows 
 {\color{black}us} to 
compute an approximate GBNE.
Specifically, we adopt the well-known 
polynomial decision rule 
 (see e.g.~\cite{bampou2012polynomial})
to approximate the GBNE
and subsequently convert the GBNE model into 
the stochastic generalized Nash equilibrium (SGNE) model,
 {\color{black} which can be solved by a sample average approximation scheme; see \cite{Franci2022Stochastic, guo2021existence,ravat2017existence}.}
This  {\color{black}polynomial approximation} approach has been proposed by Guo et al.~\cite{guo2021existence}
for classic BNE models.
When we use  {\color{black}this} approach to tackle 
the GBNE model, we need to address  
new challenges arising from 
type/action dependent 
feasible sets.
To ease the exposition,
we  confine 
our discussion to
the case where 
 {\color{black}the}
type space 
$\Theta_i\subset\R$ is a convex and compact set
and 
 {\color{black}
the action space $\mathbb{A}_i=[\underline{\underline{a}}_i, \bar{\bar{a}}_i]\subset\R$, where $\underline{\underline{a}}_i, \bar{\bar{a}}_i$ are constants for $i\in N$.
}

\subsection{Feasibility and existence of 
polynomial-GBNE}
We begin by  restricting player
$i$'s response function $f_i(\theta_i)$ to a polynomial function, written 
$\sum_{l=0}^d v_{i,l} \theta_i^l$, for $i\in N$.
 {\color{black} Here $\theta_i\in \Theta_i \subset \R$ and $\theta_i^l$ means $\theta_i$ to the power of $l$.
Let 
$
\xi_d(\theta_i):=(1,\theta_i,\theta_i^2,\dots,\theta_i^d)^T
$ 
be a vector of monomials in $\theta_i$ to the power from $l=0$ to $l=d$ at fixed point $\theta_i$.}
Then we can write 
$f_i(\theta_i)$ as 
$v_i^T\xi_d(\theta_i)$,  where $v_i=(v_{i,0},\dots,v_{i,d})$.
In this case, $f_i(\theta_i)$ is uniquely determined by $v_i$.
To facilitate the convergence analysis
in the forthcoming discussions, we 
restrict 
the polynomial response functions to the case that $v_i\in \mathbb{V}_i^d$
and
\bgeq
\mathbb{V}_i^d := \left\{ v_i\in\R^{d+1}: \underline{\underline{a}}_i\leq v_i^T \xi_d(\theta_i) \leq \bar{\bar{a}}_i, \forall \theta_i \in \Theta_i\right\},
\edeq
 {\color{black}where we assume the action space as $\mathbb{A}_i = [\underline{\underline{a}}_i,\bar{\bar{a}}_i]\subset \R$ as discussed above.} 
Let  $\mathbb{V}^d:= \mathbb{V}_1^d\times\cdots\times \mathbb{V}_n^d$,
\bgeqn
S_i^d:=\{v_i^T\xi_d(\theta_i): v_i\in \mathbb{V}_i^d \subset \R^{d+1}\}
\; \inmat{and} \; S^d:= S_1^d\times\cdots\times S_n^d.
\edeqn

\begin{definition}
\label{def:poly_GBNE}
    A polynomial-GBNE is an $n$-tuple $f^d:=(f^d_1,\dots,f^d_n)\in S^d$ mapping from $\Theta_1\times\cdots\times\Theta_n$ to $\mathbb{A}_1\times\cdots\times\mathbb{A}_n$ such that 
    \bgeqn
    \label{eqn:poly_gbne}
    \begin{aligned}
    f^d\in \arg \max_{\tilde{f}^d\in S^d}& \sum_{i=1}^n \mathbb{E}_{\eta} [u_i(\tilde{f}^d_i(\theta_i),f^d_{-i}(\theta_{-i}),\theta)]\\
    \inmat{s.t.}&  \int_{\Theta_{-i}} g_i(\tilde{f}^d_i(\theta_i),f^d_{-i}(\theta_{-i}),\theta)d\eta_i(\theta_{-i}|\theta_i)\leq 0,\forall \theta_i\in\Theta_i, i \in N , 
    \end{aligned}
    \edeqn
    where $f^d_i$ is player $i$’s best polynomial response strategy with maximum degree $d$. 
\end{definition}
Problem (\ref{eqn:poly_gbne})  can be rewritten as 
    \bgeqn
    \label{eqn:poly_gbne'}
    \begin{aligned}
        V\in \arg \max_{\tilde{V}\in \mathbb{V}^d}& \sum_{i=1}^n \mathbb{E}_{\eta} [u_i(\tilde{v}_i^T \xi_d(\theta_i),v_{-i}^T\xi_d(\theta_{-i}),\theta)]\\
        \inmat{s.t.}&  \int_{\Theta_{-i}} g_i(\tilde{v}_i^T\xi_d(\theta_i),v_{-i}^T\xi_d(\theta_{-i}),\theta)d\eta_i(\theta_{-i}|\theta_i)\leq 0, \forall \theta_i\in\Theta_i, i \in N ,
    \end{aligned}
    \edeqn
    where 
    $V = (v_1,\dots,v_n)\in \mathbb{V}^d\subset \R^{n\times(d+1)}$.
From the definition, we can see
that a polynomial-GBNE 
 {\color{black}can}
be viewed 
as an approximation of GBNE in problem (\ref{eqn:GBNE}) where each player's response function is restricted to a polynomial with degree $d$. 
 {\color{black}By} comparison 
with (\ref{eqn:GBNE}), (\ref{eqn:poly_gbne'}) is a 
finite-dimensional 
SGNE problem.
Let
\bgeqn
\label{eqn:feasible_set_of_polynomial}
\;\mathcal{V}_i^d(v_{-i},\theta_i) &:=& \left\{ \tilde{v}_i\in 
\mathbb{V}_i^d
: \int_{\Theta_{-i}} g_i(\tilde{v}_i^T\xi_d(\theta_i),v_{-i}^T\xi_d(\theta_{-i}),\theta)d\eta_i(\theta_{-i}|\theta_i)\leq 0 \right\}
\edeqn
and $\mathcal{V}_i^d(v_{-i}) := \bigcap\limits_{\theta_i\in \Theta_i} \mathcal{V}_i^d(v_{-i},\theta_i)$.
The next lemma states that the inequality 
system in (\ref{eqn:poly_gbne'}) 
satisfies Slater's condition under some moderate conditions.

\begin{lemma}
\label{lemma:slater_condition_poly}
    Assume: (a) $g_i(a_i,a_{-i},\theta), i \in N$ are H{\"o}lder continuous in {\color{black}$a_i$} and
    $\theta \in \Theta$: there are positive parameters $L_g>0$ and {\color{black}$q_{g,a},q_{g,\theta}\in (0,1]$ such that for all $a_{-i}\in \mathbb{A}_{-i}$
    \bgeq
        |g_i(a_i',a_{-i},\theta')-g_i(a_i'',a_{-i},\theta'')| \leq L_g (\| a_i' - a_i'' \|^{q_{g,a}} + \| \theta'-\theta''\|^{q_{g,\theta}}),&\\ 
        \forall a_i', a_i'' \in \mathbb{A}_i, \theta',\theta''\in \Theta_i;&
    \edeq
    }
  (b)  for $i\in N$, $g_i(a_i,a_{-i},\theta)$ is convex in $a_i\in \mathbb{A}_i$;
    (c)  {\color{black}the} inequality system in 
    model (\ref{eqn:GBNE})
    satisfies Slater's condition {\color{black}:} 
     {\color{black} there exists a positive constant $\alpha > 0$ such that for any $(f_{-i},\theta_i)$, there exists a feasible solution, denoted by $\bar{a}_i(f_{-i},\theta_i)\in \mathbb{A}_i$, satisfying
        $h_i(\bar{a}_i(f_{-i},\theta_i),f_{-i},\theta_i)\leq -\alpha$;
        }
        
    {\color{black}
    Then,
    for any positive constant $\alpha'\in (0,\frac{\alpha}{2})$,
    there exists a positive number $M$ and a sequence of polynomials $\{\bar{f}^d\}_{d\geq M}$ such that for any $f_{-i}\in \mathcal{C}_{-i}$ and $d\geq M$,
    \bgeqn
    \label{eqn:slater_poly}
        h_i(\bar{f}_i^d(\theta_i),f_{-i},\theta_i)\leq -\alpha', \forall \theta_i\in \Theta_i.
    \edeqn
    Moreover, for any $v_{-i}\in \mathbb{V}_{-i}^d$ and $f_{-i}:= v_{-i}^T\xi_d(\cdot)$,
    (\ref{eqn:slater_poly}) still holds.
    }
    
\end{lemma}
\textbf{Proof.}
Since the type space $\Theta_i$ is 
a convex and compact interval, 
we assume with loss of generality 
that $\Theta_i = [0,1]$.
{\color{black} For any fixed $f_{-i}\in \mathcal{C}_{-i}$}, under condition (c), it follows by Lemma \ref{lemma:hoffman} (iii) that there exists a continuous function $\bar{f}_i:\Theta_i\rightarrow \mathbb{A}_i$ and a positive constant $\alpha$ such that for $i\in N$,
$
h_i(\bar{f}_i(\theta_i),f_{-i},\theta_i)\leq  -\frac{\alpha}{2},\; \forall \theta_i\in\Theta_i.
$
By the Weierstrass theorem, we can find a sequence of Bernstein polynomials \linebreak $\{B^d(\theta_i; \bar{f}_i)\}_{d=1}^\infty$ 
approximating $\bar{f}_i$, where
\bgeq
B_i^d(\theta_i; \bar{f}_i) := \sum_{j=1}^d \bar{f}_i(j/d) \tbinom{d}{j} \theta_i^j(1-\theta)^{d-j}, d=1,2,\cdots.
\edeq
For $\alpha'\in (0,\frac{\alpha}{2})$, let $\epsilon := \frac{\alpha}{2}-\alpha'$. Then there exists a 
positive number 
$M$ sufficiently large 
such that for all $d\geq M$, 
$\|B_i^d(\cdot; \bar{f}_i) - \bar{f}_i\|_\infty \leq \left(\frac{\epsilon}{L_g}\right)^{1/q_{q,a}}$, for $ i\in N$.
By the H\"older continuity of $g_i(a_i,a_{-i},\theta), i \in N$ in $a_i$, we have
\bgeq
&& \int_{\Theta_{-i}}g_i(B_i^d(\theta_i; \bar{f}_i), 
f_{-i}(\theta_{-i}),
\theta)  {\color{black}d\eta_i}(\theta_{-i}|\theta_i) \\
&\leq&  \int_{\Theta_{-i}} g_i(\bar{f}_i(\theta_i), 
f_{-i}(\theta_{-i}),
\theta)  {\color{black}d\eta_i} (\theta_{-i}|\theta_i) + L_g \| B_i^d(\cdot; \bar{f}_i) - \bar{f}_i \|^{q_{g,a}}\\
&\leq & -\frac{\alpha}{2} + L_g \left(\frac{\epsilon}{L_g}\right)^{q_{g,a}/q_{g,a}}  = - \alpha'.
\edeq
Let $\bar{f}_i^d = B_i^d(\cdot; \bar{f}_i), 
i\in N$
and $d\geq M$.
{\color{black} 
The conclusion follows by 
setting $f_{-i}:= v_{-i}^T\xi_d(\cdot)$
where $v_{-i}\in \mathbb{V}_{-i}^d$ because  
$v_{-i}^T\xi_d(\cdot)\in \mathcal{C}_{-i}$}.
\hfill $\Box$

Now, we are ready to state  {\color{black}the} existence of polynomial GBNE.

\begin{theorem}
    Assume the setting and conditions of Lemma \ref{lemma:slater_condition_poly}. Assume additionally: (a)
    $\Theta_i$ is  convex and compact;
        (b) $u_i(a_i,a_{-i},\theta)$ is continuous in $(a_i,a_{-i},\theta)$;
        (c) for each 
        $(a_{-i},\theta)$, $u_i(a_i,a_{-i},\theta)$ is concave on $\mathcal{A}_{i}(f_{-i},\theta_i)$.
Then  problem (\ref{eqn:poly_gbne'}) has a polynomial GBNE.
\end{theorem}

\noindent
\textbf{Proof.} Write the objective function of 
program (\ref{eqn:poly_gbne'}) concisely 
as 
\bgeq
\Gamma(\tilde{V},V) := 
    \sum_{i=1}^n \mathbb{E}_{\eta} [u_i(\tilde{v}_i^T \xi_d(\theta_i),v_{-i}^T\xi_d(\theta_{-i}),\theta)],
\edeq
where $V :=(v_1,\dots,v_n)$.
Since $u_i(a_i,a_{-i},\theta)$ is continuous in $(a_i,a_{-i},\theta)$ and concave in $a_i$, then $\Gamma(\tilde{V},V)$ is continuous and concave in $\tilde{V}$ for fixed $V$.
Let $\mathcal{V}_i^d(v_{-i})$ be defined below (\ref{eqn:feasible_set_of_polynomial}).
Consider the set-valued mapping 
\bgeq
\Phi(V) := \left\{ \tilde{V}\in \mathbb{V}^d :
\tilde{V} = \arg \max\limits_{\tilde{V}\in \mathbb{V}^d} \Gamma(\tilde{V},V) \text{ s.t. } \tilde{v}_i\in 
\mathcal{V}_i^d(v_{-i}), \forall i\in N\right\}.
\edeq
Observe that the polynomial-GBNE is a fixed point of $\Phi(\cdot)$. Next, we 
prove that the fixed point of $\Phi(\cdot)$ exists.

Since $g_i(a_i,a_{-i},\theta)$ is continuous in $(a_i,a_{-i},\theta)$ and is convex in $a_i$, then 
$\mathcal{V}_i^d(v_{-i},\theta_i)$ is compact for all $\theta_i$.
Thus, the feasible set $\mathcal{V}_i^d(v_{-i})= \bigcap\limits_{\theta_i\in \Theta_i} \mathcal{V}_i^d(v_{-i},\theta_i)$ is non-empty, convex and compact. 
On the other hand, by Lemma \ref{lemma:slater_condition_poly}, Slater's condition for the inequality system in polynomial-GBNE model (\ref{eqn:poly_gbne'}) is satisfied.
By virtue of Lemma~\ref{lemma:hoffman}(ii) (where 
the correspondence of the two parametric inequality systems 
in (\ref{eqn:para_min_problem}) and (\ref{eqn:poly_gbne'}) can be identified as 
$\tilde{v}_i \leftrightarrow x$, $v_{-i}\leftrightarrow t$ and $1 \leftrightarrow y$),
$\mathcal{V}_i^d(v_{-i})$ is continuous in $v_{-i}$. 
By 
Theorem 4.2.1 in \cite{bank1984Nonlinear},
the optimal solution set $\Phi(V)$ is u.s.c. in $V$.
Moreover, the optimal solution set $\Phi(V)$ is non-empty, convex and closed by the concavity of $\Gamma(\tilde{V},V)$ in $\tilde{V}$ and convexity and compactness of $\mathcal{V}_i^d(v_{-i}), i\in N$. 
Therefore, the mapping $\Phi(V)$ satisfies all  {\color{black}the} conditions for Kakutani's fixed-point theorem (see e.g. Theorem 3.2.2 in \cite{ichiishi1983game}) and thus there exists a $V^*$ such that $V^* \in \Phi(V^*)$. The proof is complete.
\hfill $\Box$

\subsection{Convergence of polynomial-GBNEs}
In this section, we discuss
the convergence of  {\color{black}the} polynomial-GBNEs obtained 
 {\color{black}by} solving (\ref{eqn:poly_gbne}) as the degree $d$ increases.
The next theorem addresses this.

\begin{theorem}
\label{thm:convergence_poly}
    Let 
    $\{f^d\}_{d=1}^\infty$ be a sequence of 
    polynomial-GBNEs obtained 
    from solving 
    problem (\ref{eqn:poly_gbne}). Assume:
    (a) for $i\in N$,
    $u_i(a_i,a_{-i},\theta)$ are 
    continuous 
    in $(a_i,a_{-i},\theta)$, 
    and
    {\color{black}
    $g_i(a_i,a_{-i},\theta), i\in N$ are H{\"o}lder continuous in $a_i$ and $\theta$ with constants $L_g>0,q_{g,a}\in (0,1]$, $q_{g,\theta}\in (0,1]$;
    }
    (b)  {\color{black}for $i\in N$,} $g_i(a_i,a_{-i},\theta)$ is convex in $a_i$;
    (c) for any fixed $f_{-i}\in \mathcal{F}_{-i}$, Slater's condition holds for the parametric inequality system of the GBNE model (\ref{eqn:GBNE}) with 
    a constant $\alpha>0$;
    (d) $\mathbb{A}_i, i\in N$ are compact with $\inmat{diam}(\mathbb{A}_i) \leq A, \forall i\in N$.
    Then the following assertions hold.
    \begin{enumerate}
        \item[(i)] For 
        any continuous strategy function $f\in \mathcal{C}$,
        any feasible response function $\tilde{f}\in \mathcal{C}(f)$ under definition (\ref{def:cont_func_set})  
        and any constant $\epsilon\in (0,1]$, there exist a positive number $M$ sufficiently large and a sequence of polynomials $\{\hat{f}^d\}_{d\geq M}$
        {\color{black} feasible to GBNE model (\ref{eqn:GBNE}), e.g., $\{\hat{f}^d\}_{d\geq M}\subset \mathcal{C}(f)$} 
        such that for $i\in N$,
        \bgeq
            \| \hat{f}_i^d-\tilde{f}_i \|_\infty \leq  \left(\frac{C}{\alpha^\prime} L_g+1\right) \epsilon^{q_{g,a}},
        \edeq
        where $C=\inmat{diam}(\cup_{\theta_i\in\Theta_i}\mathcal{A}_i(f_{-i},\theta_i))+ 1$  and $\alpha'\in (0,\frac{\alpha}{2})$.
        \item[(ii)] Every cluster point of $\{f^d\}_{d=1}^\infty$, written $f^\infty$,  is a GBNE of model (\ref{eqn:GBNE''2}).
    \end{enumerate}


\end{theorem}


\textbf{Proof.}
Part (i).
Since the type space $\Theta_i$ is 
a convex and compact interval, 
we assume with loss of generality 
that $\Theta_i = [0,1]$. 
For any fixed 
$f\in \mathcal{C}$ and feasible 
$\tilde{f}_i\in \mathcal{C}_i(f_{-i}), i\in N$, 
we construct a sequence of Bernstein polynomials $\{B^d(\theta_i; \tilde{f}_i)\}_{d=1}^\infty$ approximating $\tilde{f}_i$ 
$B_i^d(\theta_i; \tilde{f}_i) := \sum_{j=1}^d \tilde{f}_i(j/d) \tbinom{d}{j} \theta_i^j(1-\theta)^{d-j}$.
By the Weierstrass theorem, for any $\epsilon\in (0,1]$, there is a constant $M'>0$ such that for all $d\geq M'$, 
$
\|B_i^d(\theta_i; \tilde{f}_i) - \tilde{f}_i\|_\infty \leq \epsilon$.
By the H\"older continuity of $g_i(\cdot,a_{-i},\theta)$, 
\bgeq
&&\int_{\Theta_{-i}} g_i(B_i^d(\theta_i; \tilde{f}_i),
f_{-i}(\theta_{-i}),
\theta) d\eta_i(\theta_{-i}|\theta_i)\\
\leq&& \int_{\Theta_{-i}} g_i(\tilde{f}_i(\theta_i),
f_{-i}(\theta_{-i}),
\theta) d\eta_i(\theta_{-i}|\theta_i)+L_g \| B_i^d(\cdot; \tilde{f}_i) - \tilde{f}_i \|_\infty^{q_{g,a}} \leq L_g \epsilon^{q_{g,a}}, \quad \forall \theta_i\in \Theta_i,
\edeq
where the second inequality is 
 {\color{black}because} $\tilde{f}$ is 
feasible 
 {\color{black}for}
the GBNE model (\ref{eqn:GBNE}). 
On the other hand, by Lemma~\ref{lemma:slater_condition_poly}, for any $\alpha'\in (0,\frac{\alpha}{2})$, there exists a positive constant $M''$ and a sequence of $\{\bar{f}^d\}_{d\geq M''}$
such that 
for all $d\geq M''$ and $i\in N$,
\bgeq
\int_{\Theta_{-i}} g_i(\bar{f}_i^d(\theta_i),
f_{-i}(\theta_{-i})
,\theta) d\eta_i(\theta_{-i}|\theta_i)\leq  -\alpha',\; \forall \theta_i\in\Theta_i.
\edeq
Let 
$
\gamma  :=  \max\limits_{i\in N} \left\{ \sup\limits_{\theta_i\in \Theta_i}\left(\int_{\Theta_{-i}} g_i(B_i^d(\theta_i; \tilde{f}_i),
f_{-i}(\theta_{-i})
,\theta) d\eta(\theta_{-i}|\theta_i)\right)_+\right\}
 \in  [0,
L_g\epsilon^{q_{g,a}}]
$
and 
$\hat{f}_i^d := (1-\frac{\gamma}{\gamma+\alpha'})B_i^d(\cdot;\tilde{f}_i) + \frac{\gamma}{\gamma+\alpha'} \bar{f}_i^d$.
By virtue of 
Lemma~\ref{lemma:hoffman} (i),
for any $d\geq M := \max\{M',M''\}$, $\hat{f}_i^d$ is feasible 
 {\color{black}for} 
the GBNE model (\ref{eqn:GBNE}) and
\bgeq
\| \hat{f}_i^d-B_i^d(\cdot;\tilde{f}_i) \|_\infty &\leq & \frac{\gamma}{\alpha^\prime}\cdot\|B_i^d(\cdot;\tilde{f}_i)-\bar{f}_i^d \|_\infty \\
&\leq & \frac{\gamma}{\alpha^\prime}\cdot
\left(\|B_i^d(\cdot;\tilde{f}_i)- \tilde{f}_i\|_\infty+\|\tilde{f}_i-\bar{f}_i^d \|_\infty\right)\\
&\leq & \frac{\gamma}{\alpha^\prime}\cdot\left[\epsilon + \inmat{diam}(\cup_{\theta_i\in\Theta_i}\mathcal{A}_i(f_{-i},\theta_i))\right] 
\leq \frac{C\gamma}{\alpha^\prime},  
\edeq
where 
$C := 1 + \inmat{diam}(\cup_{\theta_i\in\Theta_i}\mathcal{A}_i(f_{-i},\theta_i))$.
Therefore, 
\bgeq
\| \hat{f}_i^d-\tilde{f}_i \|_\infty & \leq & \| \hat{f}_i^d-B_i^d(\cdot;\tilde{f}_i) \|_\infty + \| B_i^d(\cdot;\tilde{f}_i)^d-\tilde{f}_i \|_\infty \nonumber\\
& = & \frac{C}{\alpha'}\gamma + \epsilon \leq \left(\frac{C}{\alpha^\prime} L_g+1\right) \epsilon^{q_{g,a}}, \quad \forall i \in N.
\edeq

{\color{black}
Part (ii). By taking a subsequence if necessary, we assume for the simplicity of exposition that  
$\|f^d-f^\infty\|_\infty \rightarrow 0$ as $d\rightarrow \infty$. 
By the definition of the polynomial-GBNE $f_i^d$ and the continuity of $g_i(a_i,a_{-i},\theta)$, we have
$\int_{\Theta_{-i}} g_i(f_i^d(\theta_i), 
f_{-i}^d(\theta_{-i}),\theta) d\eta_i(\theta_{-i}|\theta_i) \leq 0$,
and thus
$
\int_{\Theta_{-i}} g_i(f_i^{\infty}(\theta_i),
f_{-i}^\infty(\theta_{-i}),\theta) d\eta_i(\theta_{-i}|\theta_i) \leq 0,
$
which means $f^\infty$ is feasible 
 {\color{black}for}  {\color{black}the} 
GBNE model (\ref{eqn:GBNE''2}).

\begin{figure}[t]
    \centering
    \includegraphics[width=0.3\textwidth]{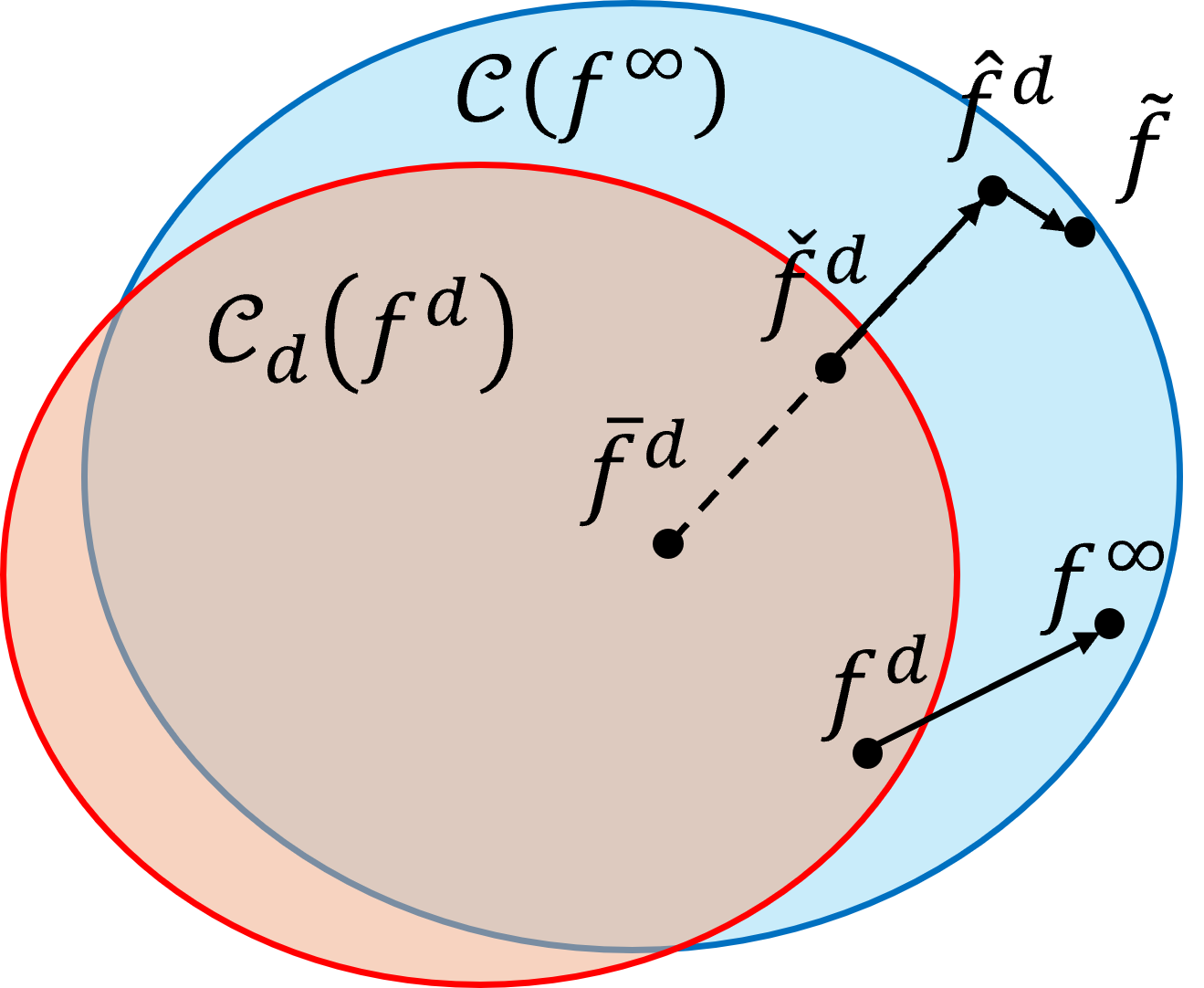}
    \caption{
    \small Illustration the proof of Part (ii), where $\mathcal{C}_d(f^d)$ denotes the polynomial feasible set defined by  $f^d$.}
    \label{fig:polynomial}
    \vspace{-15pt}
\end{figure}

For any 
$\tilde{f} \in \mathcal{C}(f^\infty)$, 
it follows from Part (i) that
there exists
a sequence of polynomial functions
$\{\hat{f}^d\}_{d\geq M}\subset \mathcal{C}(f^\infty)$ 
 {\color{black}that converge}
to $\tilde{f}$.
Let 
\bgeq
\gamma_i^d := \left(\int_{\Theta_{-i}} g_i(\hat{f}_i^d,
f_i^d(\theta_{-i})
,\theta) d\eta(\theta_{-i}|\theta_i)\right)_+.
\edeq
Then
$
\gamma_i^d \rightarrow \left(\int_{\Theta_{-i}} g_i(\tilde{f}_i(\theta_i),
f_{-i}^\infty(\theta_{-i})
,\theta) d\eta_i(\theta_{-i}|\theta_i)\right)_+ = 0
$
as $d\rightarrow +\infty$.
On the other hand, by Lemma \ref{lemma:slater_condition_poly}, for each $f^d$ with $d$
being large enough,
there is a positive parameter $\alpha'$ and a polynomial $\bar{f}^d$ with degree $d$ such that 
\bgeq
\int_{\Theta_{-i}} g_i(\bar{f}_i^d(\theta_i),
f_{-i}^d(\theta_{-i})
,\theta) d\eta_i(\theta_{-i}|\theta_i)\leq  -\alpha',\; \forall \theta_i\in\Theta_i.
\edeq
Let 
$\check{f}_i^d := (1-\frac{\gamma_i^d}{\gamma_i^d+\alpha'})\hat{f}_i^d + \frac{\gamma_i^d}{\gamma_i^d+\alpha'} \bar{f}_i^d$.
 {\color{black}Like}
the proof of Lemma~\ref{lemma:hoffman}, we have
\bgeq
\int_{\Theta_{-i}} g_i(\check{f}_i^d(\theta_i),
f_{-i}^d(\theta_{-i})
,\theta) d\eta_i(\theta_{-i}|\theta_i)\leq  0,\; \forall \theta_i\in\Theta_i,
\edeq
which means  {\color{black}that} $\check{f}^d$ lies 
 {\color{black}within}
the feasible set of polynomial response functions parameterized  by $f^d$ (denoted by  ${\cal C}_d(f^d)$),
and 
$
\|\hat{f}_i^d-\check{f}_i^d\|_\infty
\leq
\frac{A}{\alpha'}\gamma_i^d.
$
Therefore, 
$
\check{f}_i^d \rightarrow \hat{f}_i^d \rightarrow \tilde{f}_i
$
as $d\rightarrow \infty$.
By the definition of polynomial-GBNEs 
$\{f^d\}_{d=1}^\infty$ in (\ref{eqn:poly_gbne}),
\bgeq
\mathbb{E}_{\eta} [ u_i(\check{f}_i^d(\theta_i),f_{-i}^d(\theta_{-i}),\theta_i) ] \leq \mathbb{E}_{\eta} [ u_i(f_i^d(\theta_i),f_{-i}^d(\theta_{-i}),\theta_i) ].
\edeq
By driving $d$ to infinity and 
the continuity of $u_i(a_i,a_{-i},\theta)$,
we obtain
\bgeq
\mathbb{E}_{\eta} [ u_i(\tilde{f}_i(\theta_i),f_{-i}^\infty(\theta_{-i}),\theta_i) ] \leq \mathbb{E}_{\eta} [ u_i(f_i^\infty(\theta_i),f_{-i}^\infty(\theta_{-i}),\theta_i) ],
\edeq
which implies that $f^\infty$ is a GBNE.
\hfill $\Box$
}

Theorem \ref{thm:convergence_poly} shows that 
any cluster points of
the sequence of polynomial-GBNEs $\{f^d\}_{d}$ are 
GBNEs of model (\ref{eqn:GBNE}).
The theorem provides theoretical 
 {\color{black}grounds}
for us to obtain an approximate GBNE
by solving the polynomial-GBNE 
problem (\ref{eqn:poly_gbne}). 

\section{Numerical tests}
    \label{sec:numerical}
 {\color{black} In Section \ref{sec:polynomial},
a polynomial decision rule 
 {\color{black}was}
employed to approximate the GBNE model (\ref{eqn:GBNE}) with the SGNE model (\ref{eqn:poly_gbne'}). 
The polynomial decision rule allows us to restrict the feasible set from the set of continuous functions to the set of $d$-degree polynomial functions (which is uniquely determined by $d+1$ coefficients).
}
To 
investigate the performance of 
the GBNE model and the polynomial approximation scheme, 
we 
have conducted some numerical experiments 
by applying them to Example \ref{exam:cournot_game}.
 {\color{black}
Note that model (\ref{eqn:poly_gbne'}) contains
    mathematical expectations in both the objective and the constraints 
    which could be difficult to calculate either because the true probability distribution of $\eta$ is unknown
in some data-driven problems, or because
computing the integrals could be prohibitively expensive.
    The well-known sample average approximation (SAA) method 
is subsequently proposed  to solve the SGNE problem (\ref{eqn:poly_gbne'}) as follows:
    \begin{eqnarray}
        \label{eqn:poly_gbne_saa}
        \begin{aligned}
        V\in \arg \max_{\tilde{V}\in \mathbb{V}^d}& \sum_{i=1}^n \frac{1}{M} \sum_{j=1}^M [u_i(\tilde{v}_i^T \xi_d(\theta_i^j),v_{-i}^T\xi_d(\theta_{-i}^j),\theta^j)]\\
        \text{s.t.}&  \frac{1}{M}  \sum_{j=1}^M g_i(\tilde{v}_i^T\xi_d(\theta_i^j),v_{-i}^T\xi_d(\theta_{-i}^j),\theta^j)d\leq 0, i \in N,
        \end{aligned}
    \end{eqnarray}
where $\{\theta^j\}_{j=1}^M$ are iid samples.
Convergence of SAA (applied to stochastic games/ equilibrium problems) is well documented, see \cite{ravat2011characterization,xu2010uniform,xu2013stochastic} for
stochastic equilibrium problems
and \cite{guo2021existence} 
for classic BNE models.
The experiments are carried out using the
standard software GAMS with an extended mathematical programming framework devised by Kim and Ferris \cite{kim2019solving}. The general algebraic modeling system (GAMS) is a powerful software and has been used widely to solve 
optimization and game theoretic problems.
}


\begin{example} [Symmetric  {\color{black}two-player} Bayesian Cournot game]
\label{example:2p_sym_BCG}
    For each player $i=1,2$, 
    action $a_i$ is the quantity of products to be produced, and the type $\theta_i$ is the quantity of resources required for one unit of the product, 
    where 
    the type spaces are $\Theta_1=\Theta_2=[\underline{\theta},\overline{\theta}] = [1, 2]$.
    The cost function 
    is $c_i(a_i, \theta_i) = k_i \theta_i a_i$, w.l.o.g. $k_i=1$.
    The inverse demand function is defined as $p(a_1,a_2) = A - a_1- a_2$, where $A=100$.
    The type $\theta = (\theta_1,\theta_2)$ is drawn from a uniform distribution $\eta$ supported over $\Theta_1\times \Theta_2 = [1,2]\times [1,2]$.  
    The  {\color{black}upper bound of the} feasible resource for each player is 
    $R_i=30,i=1,2$ 
    respectively:
    $
    \theta_i a_i \leq R_i, i=1,2
    $.
    The coupled constraints 
    mean that  {\color{black}the upper bound of} the level of the 
    shared resource available to the industry is 
    $R_c=50$: 
    \bgeqn
    \label{cons:coupled_numerical}
    \theta_1 a_1 + \theta_2 a_2 \leq R_c.
    \edeqn
    We approximate the GBNE model by polynomial functions with degree $d$ and the sample average approximation with the sample size $M=100$ for each type space $\Theta_i$.
    The empirical distribution $\eta^M$ is supported by $\Theta_1^M\times\Theta_2^M := \{(\theta_1^{j_1}, \theta_2^{j_2})\}_{j_1,j_2=1}^M$, with marginal distributions  $\eta^{M}_{-1}$ over $\Theta_2^M := \{\theta_2^{j_2}\}_{j_2=1}^M$, and $\eta^{M}_{-2}$ over $\Theta_1^M := \{\theta_1^{j_1}\}_{j_1=1}^M$. 
    We 
    compare two models in this example: 
    Bayesian Cournot game with the expected coupled constraint (BCG-ECC)  {\color{black} to find a $2$-tuple $(v_1,v_2)$ such that, } 
    \bgeqn
    \label{example:BCG_ECC}
    \begin{aligned}
        v_i \in \arg \max_{v_i\in \mathcal{V}_i^d} &\; \mathbb{E}_{\eta^M} \left( A - v_1^T \xi_d(\theta_1) - v_2^T \xi_d(\theta_2) - \theta_i  \right) v_i^T \xi_d(\theta_i),\; i=1,2,\\
        \inmat{s.t.} &\; \overbrace{\mathbb{E}_{\eta^{M}_{-i}} \left( \theta_1 v_1^T \xi_d(\theta_1) + \theta_2 v_2^T \xi_d(\theta_2) - R_c \right) \leq 0}^\text{Expected Coupled Constraint},\forall \theta_i \in \Theta_i^{M}\\
        &\; \theta_i v_i^T \xi_d(\theta_i) \leq R_i, \forall \theta_i \in \Theta_i^{M},
    \end{aligned}
    \edeqn
    and 
    Bayesian Cournot game with CVaR coupled constraint (BCG-CCC)  {\color{black} to find a $2$-tuple $(v_1,v_2)$ such that, } 
    \bgeqn
    \label{example:BCG_CCC}
    \begin{aligned}
        v_i \in \arg \max_{v_i\in \mathcal{V}_i^d} &\; \mathbb{E}_{\eta^M} \left( A - v_1^T \xi_d(\theta_1) - v_2^T \xi_d(\theta_2) - \theta_i  \right) v_i^T \xi_d(\theta_i),\; i=1,2,\\
        \inmat{s.t.} &\; \overbrace{t_i(\theta_i) + \frac{1}{\beta}\mathbb{E}_{\eta_{-i}^M} \bigg( \theta_1 v_1^T \xi_d(\theta_1) +\theta_2  v_2^T \xi_d(\theta_2) - R_c - t_i(\theta_i) \bigg)_+ \leq 0}^\text{CVaR Coupled Constraint}, \forall \theta_i \in \Theta_i^{M}\\
        &\; \theta_i v_i^T \xi_d(\theta_i) \leq R_i, \forall \theta_i \in \Theta_i^{M},
    \end{aligned}
    \edeqn
    where $\beta\in (0,1)$.
\end{example}

\begin{figure}[!ht]
    \centering
    \includegraphics[width=0.9\textwidth]{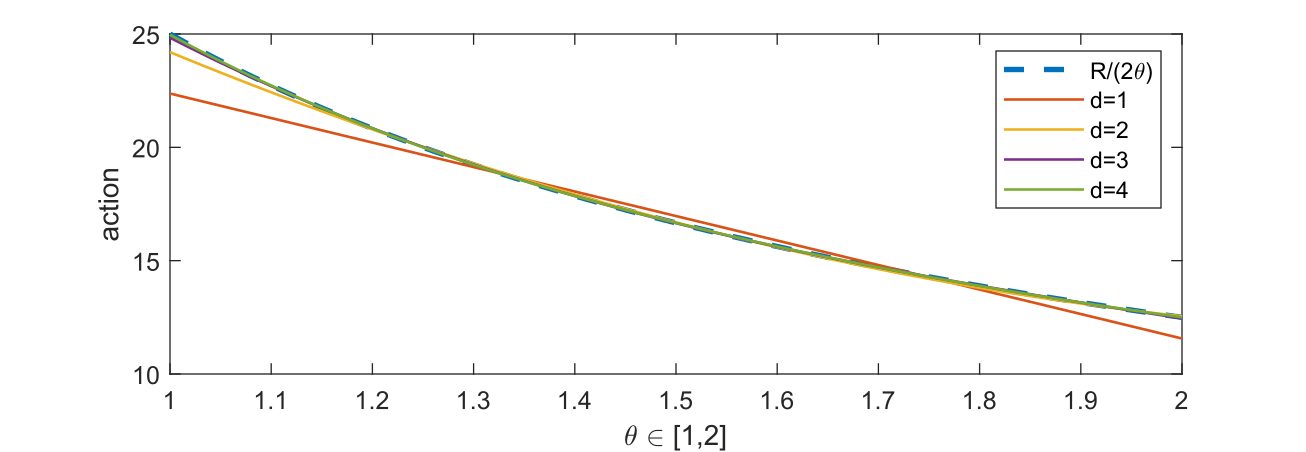}
    \caption{
    \small Optimal response functions at 
    polynomial GBNE of
    symmetric BCG with ECC.
    The blue dashed curve represents the optimal response function at the true GBNE. }
    \label{fig:ECC_d_varying}
\end{figure}

Without resource constraints and coupled constraint, the Bayesian Nash equilibrium 
of the game is
$
f_i(\theta_i) = \frac{A}{3}+\frac{\overline{\theta}-\underline{\theta}}{12}-\frac{\theta_i}{2}\approx 33.41-\frac{\theta_i}{2}, i = 1,2.
$
However, this equilibrium 
cannot be obtained 
with  {\color{black}the} constraints in (\ref{example:BCG_ECC}) or (\ref{example:BCG_CCC}). Indeed, with these constraints in place,  {\color{black}the} 
resulting GBNE is 
\bgeqn
\label{eqn:equil_BCG_with_cons}
f_i(\theta_i) = \frac{R_c}{2\theta_i} = \frac{25}{\theta_i}, \theta_i\in [1,2], i=1,2.
\edeqn
Figure \ref{fig:ECC_d_varying} visualizes the polynomial equilibria in the BCG-ECC model (\ref{example:BCG_ECC}). 
As expected, the sequence of the polynomial equilibria 
converges to the GBNE 
in (\ref{eqn:equil_BCG_with_cons}) as  {\color{black}the} degree $d$ increases.
However, the polynomial response function
at the
equilibrium in the BCG-ECC model is  {\color{black}greater} than the optimal response function at the GBNE in (\ref{eqn:equil_BCG_with_cons})
for many $\theta_i\in \Theta_i^M$, which 
implies that the coupled constraints (\ref{cons:coupled_numerical}) 
may not be satisfied with a  {\color{black}high} probability. 
For example, if both players choose a polynomial 
response function in 
the equilibrium with degree $d=1$ (red line in figure \ref{fig:ECC_d_varying}), then for $\theta_1 = 1.5$, the coupled constraint will not be satisfied with a probability of $0.58$ w.r.t. $\theta_2\in \Theta_2^M$:
$\mathbb{P}_{\theta_2\in \Theta_2^M}
(\theta_1 a_1(\theta_1) + \theta_2 a_2(\theta_2)>50|\theta_1=1.5) = 0.58$.

\begin{figure} [!ht]
     \centering
     \begin{subfigure}[b]{0.45\textwidth}
         \centering
         \includegraphics[width=\textwidth]{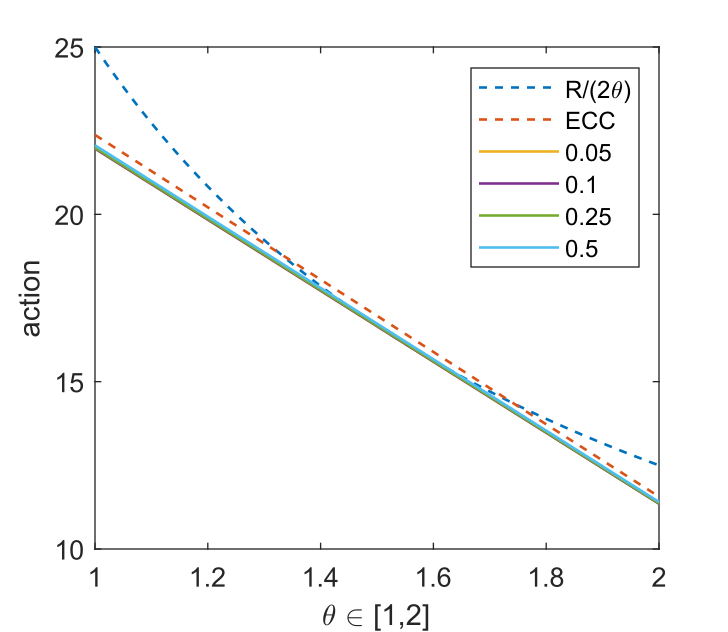}
         \caption{\small Optimal response functions at GBNE}
     \end{subfigure}
     \hfill
     \begin{subfigure}[b]{0.45\textwidth}
         \centering
         \includegraphics[width=\textwidth]{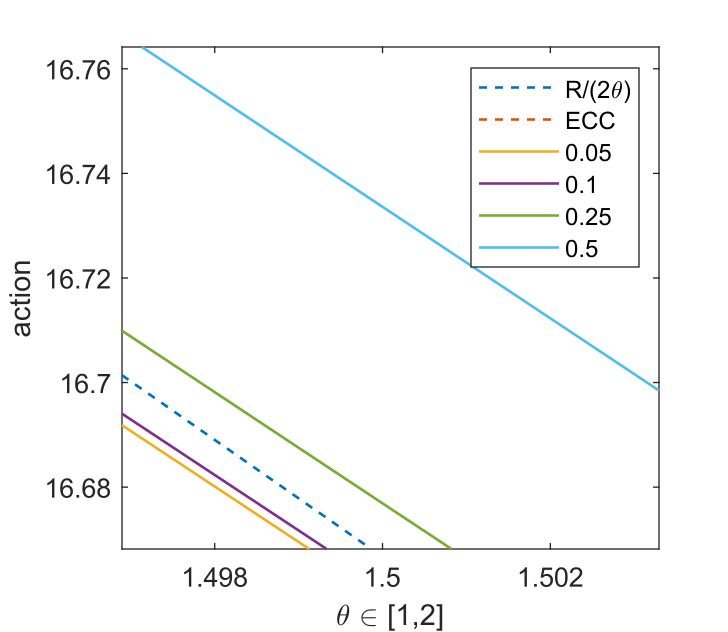}
         \caption{\small Zoom in for $\theta_i$ near $1.5$}
     \end{subfigure}
     \vspace{-0.5cm}
        \caption{\small Optimal response functions at the polynomial equilibria with degree $d=1$ for symmetric BCG-CCC model. 
    The blue dashed curve represents the optimal response function at the true GBNE. }
        \label{fig:CCC_beta_varying}
\end{figure}

To avoid  {\color{black}this} situation,
we replace 
the expected coupled constraints in the BCG-ECC model (\ref{example:BCG_ECC}) with CVaR coupled  constraints in the BCG-CCC model (\ref{example:BCG_CCC}).
To facilitate
explanation,
we set $d=1$.
Figure \ref{fig:CCC_beta_varying} depicts
the optimal response functions at the equilibrium with $\beta$  varying from
$0.05$ to $0.1, 0.25$ and $0.5$.
We can see that the red dotted line 
(representing the optimal response function in the BCG-ECC model) lies
above the other lines (representing the optimal response functions in the BCG-CCC model), which means  {\color{black}that} each firm produces less under the 
BCG-CCC model for any type $\theta_i\in [1,2]$, in other words, 
the BCG-CCC model is more conservative.
Moreover, for the  BCG-CCC model,
Figure~\ref{fig:CCC_beta_varying} (b)
shows the optimal  response functions
shift upwards as $\beta$ increases
which means that the firms are less conservative.
Table \ref{tab:probability}
displays the 
max probability of the coupled constraint being violated:  $\varrho :=\max\limits_{\theta_1\in\Theta_1^M} \mathbb{P}_{\eta_1(\cdot|\theta_1)}(\theta_1 a_1(\theta_1) + \theta_2 a_2(\theta_2)>50|\theta_1)$.

\begin{table} [!ht]
    \centering
    \caption{\small Max probability of coupled constraint not being satisfied.}
    \begin{tabular}{cccccc}
    \hline
    \multirow{2}{*}{} & \multirow{2}{*}{BCG-ECC} & \multicolumn{4}{c}{BCG-CCC} \\ \cline{3-6} 
                      &                      & $\beta=0.05$ & $\beta=0.1$ & $\beta=0.25$ & $\beta=0.5$ \\ \hline
    $\varrho$ & \multicolumn{1}{c}{0.58} & 0.04  &  0.06 &  0.12 &  0.30   \\ \hline
    \end{tabular}
    \label{tab:probability}
\end{table}

\appendix

\bibliographystyle{unsrt}
\bibliography{literature}


\end{document}